\documentclass[12pt]{article}
\usepackage{amssymb,amsmath,amsthm,epsfig}

\usepackage{eepic}
\usepackage{epic}

\usepackage{graphicx}
\usepackage{color}
\usepackage{ifpdf}

\theoremstyle{plain}
\newtheorem{lem}{Lemma}[section]

\newtheorem{prop}[lem]{Proposition}

\newtheorem{example}{Example}

\theoremstyle{definition}

\theoremstyle{remark}
\newtheorem{rem}{Remark}[section]

\begin{document}
\title{ \large\bf Mixed Multiscale Finite Volume Methods for Elliptic Problems in Two-phase
Flow Simulations}
\author{
Lijian  Jiang\thanks{IMA, University of Minnesota, Minneapolis, MN 55455.  Email: lijiang@ima.umn.edu. Corresponding
author}
\and
Ilya D.  Mishev\thanks{ExxonMobil Upstream Research Company, Houston, TX 77252. Email: Ilya.D.Mishev@exxonmobil.com}
}

\date{}
\maketitle

\begin{center}{\bf ABSTRACT}
\end{center}\smallskip

We develop a framework for constructing mixed multiscale finite volume methods
for elliptic equations with multiple scales arising from flows in porous media.
Some of the methods developed using the framework are already known \cite{jennylt03}; others are new.
New insight is gained for the known methods and extra flexibility is provided by the new methods.
We give as an example a mixed MsFV on uniform mesh in 2-D. This method uses novel multiscale velocity
basis functions that are suited for using global information, which is often needed to improve the accuracy
of the multiscale simulations in the case of continuum  scales with strong non-local features.
The method efficiently captures the small effects on a coarse grid.
We analyze the new mixed MsFV and   apply it to solve two-phase flow equations in heterogeneous  porous media.
Numerical examples demonstrate the accuracy and efficiency of the proposed method for modeling the
flows in porous media with non-separable and separable scales.

\section{Introduction}\label{introd}
Subsurface flows are often affected by heterogeneities in a wide range
 of length scales. This causes significant challenges for
subsurface flow modeling. Geological characterizations that capture
these effects are typically developed at scales that are too fine
for direct flow simulations.  Usually, upscaled or multiscale
models  are employed for such systems.
In  upscaling methods, the original model is coarsened by numerically
homogenizing  parameters (e.g., permeability). The simulation is performed using the
coarsened model, which may differ from the underlying fine-scale model.
In multiscale methods, the fine-scale information is
carried throughout the simulation and
the coarse-scale equations are generally not expressed
analytically, but rather formed and solved numerically.

Various numerical multiscale approaches for flows in porous media
have been developed during the past decade. A multiscale finite element method (MsFEM)
was  introduced  in \cite{hw97} and takes its origin
from the pioneering work \cite{bo83}.  Its main idea is to
incorporate the small-scale information into finite element
basis functions and capture their effect on the large scales
via finite element computations. The MsFEM in
\cite{hw97} shares some similarities with
a number of multiscale numerical methods, such as
residual free bubbles \cite{brezzi00},
variational multiscale method \cite{hughes98},
two-scale conservative subgrid approaches \cite{arbogast02}, heterogeneous multiscale method \cite{ee03}
and multiscale discontinuous Galerkin method \cite{sg08}.
Chen and Hou  have applied the MsFEM idea in combination with a mixed finite element  formulation
to propose a mixed MsFEM \cite{ch03}.  Recently, Arbogast et al. \cite{apwy07} used domain decomposition approach and
variational  mixed formulation to develop  a multiscale mortar mixed MsFEM.
Jenny et al. \cite{jennylt03} have used the ideas in \cite{hw97} and finite volume framework to design
a multiscale finite volume method (MsFV).  The MsFV and its variants have proved successful in reservoir simulations.

Here we develop a  framework for constructing mixed MsFV methods, which uses  ideas from the mixed finite volume methods
\cite{mishev02, mishev03, mc06}, multi-point flux approximations (MPFA) \cite{AavBBMann1:98, EdwardsRoger:99},
and mixed MsFEM. The mixed MsFV are mass conservative methods, which is an important property of the discretizations
used in subsurface flow simulations (see \cite{dsw04} for related discussion). The important feature of the mixed finite volume methods is
the direct approximation of the velocity, that is, specially constructed discrete spaces are used to approximate
the velocity unknowns. We propose a novel way to construct multiscale velocity basis functions that are well suited
for parallel computation.   Mixed MsFEM reduces the system of coupled equations for
pressure and velocity to a system only for the pressure. However, the reduction process is computational
expensive and has some restrictions when the global  mass matrix in mixed MsFEM  is large.
In the mixed MsFV, we  compute the inverse of each local mass matrix instead of global mass matrix  and get effective coarse-scale
transmissibilities. This computation is cheap and efficient.
In the MsFV proposed in  \cite{jennylt03}, two sets of multiscale basis functions are computed: the first set of
basis functions is to approximate pressure and the second set of basis functions is required to
construct a conservative fine-scale velocity field.  Only one set of multiscale basis functions
is constructed  in the mixed MsFV and the span of the basis functions are to approximate the velocity.
Piecewise constant is used for pressure basis in the mixed MsFV.   Hence the computation for
basis functions in the mixed MsFV is less expensive than the MsFV. To the best of our knowledge, the mixed MsFV
is a new numerical multiscale method.

Boundary effect is a great issue in many multiscale methods (e.g., \cite{hw97, ch03}). When
we construct the multiscale basis functions in the mixed MsFV, we can use constant boundary condition
for basis equations and obtain local multiscale basis. We find that
the  local multiscale basis  in the mixed MsFV
 has the similar merit  to the multiscale
basis  using oversampling technique developed  in \cite{hw97} and they are able to
reduce the boundary effect greatly. The mixed MsFV using local multiscale basis works
well for most multiscale problems and particularly for the case of separable scales (e.g., periodic media).
Furthermore, we can  employ   global information for the multiscale basis functions of
 the mixed MsFV. The global information usually represents  long range features of flows and is used
to construct multiscale basis functions.
The global information is  needed in the case of strong non-separable scales and renders much better accuracy
than local multiscale approaches \cite{aej08}.

The proposed mixed MsFV in some extend inherits the advantages of the mixed MsFEM
and MsFV  and alleviates the drawbacks of them without increase of the computational
cost. Moreover, this method and its generalizations are well suited for computation on unstructured grids
and can be easily incorporated in production reservoir simulators.
For example, to construct the velocity basis no geometric information is necessary. The pressure basis
can be computed using only the local matrices if it consists of discrete harmonic functions as demonstrated in
\cite{dmz09}. Detailed description is provided in \cite{patent_MMsFV}.

The rest of the paper is organized as follows.   Section 2 is devoted to
formulating a standard mixed finite volume method for a model elliptic equation.
In Section 3, we apply the methodology of the mixed finite volume method to develop
a mixed MsFV method. Here we design a new multiscale velocity basis function,
analyze  the proposed mixed MsFV,  and address some computational issues.
In Section 4, we generalize the approach applied in Section 3 to derive our first mixed MsFV method
and show how given a mixed FV method we can derived a corresponding mixed MsFV.
Several examples are given to demonstrate how the framework can be used to develop new mixed MsFV methods.
In Section 5, we apply a mixed MsFV  to incompressible
two-phase flows in porous media with continuum scales and separable scales.
Finally, some comments and conclusions are made.

\section{Mixed finite volume method formulation for a model problem}
\label{key-MFV}

We first define notations for  function spaces used in the paper.
\begin{eqnarray*}
L^2(\Omega)&=&\{f(x) |\int_{\Omega} |f(x)|^2 dx<\infty\},\\
H^1(\Omega)&=&\{f(x) |f(x)\in L^2(\Omega)\quad  \text{and} \quad \nabla f(x)\in [L^2(\Omega)]^d\},\\
H(div, \Omega) &=&\{f(x) |f(x) \in [L^2(\Omega)]^d \quad  \text{and} \quad  \text{div} (f(x))\in L^2(\Omega)\}.
\end{eqnarray*}

We consider  the following model elliptic  equation,
\begin{eqnarray}
\label {model-eq}
\begin{cases}
\begin{split}
-\text{div}(k(x)\nabla p)&=f(x) \ \   \text {in}  \  \   \Omega,\\
k(x)\nabla p\cdot n &=0 \ \  \text{on}  \ \  \partial \Omega,\\
\int_{\Omega} p dx&=0,
\end{split}
\end{cases}
\end{eqnarray}
where $\Omega$ is a domain in $\mathbb{R}^d$, $d = 2$ or $3$ and
$f \in L^2(\Omega)$.
Equation \eqref{model-eq} is used to model many physical processes, for example, fluid flow in porous media.
The coefficient $k(x)$ represents the permeability
and is often heterogeneous. Here $p$ represents pressure.

We define velocity $u(x)=-k(x)\nabla p$.   To simplify presentation, we shall not write
the spatial variables  $x$ for  functions  when no ambiguity occurs.
Then (\ref{model-eq})  can be rewritten as a first order system
\begin{eqnarray}
\label{pu-eq}
\begin{cases}
\begin{split}
k^{-1}u+\nabla p &=0\\
\text{div}(u) &=f.
\end{split}
\end{cases}
\end{eqnarray}
The weak mixed formulation for (\ref{pu-eq}) is:\\
\indent
Find $\{u, p\} \in \mathcal{U} \times \mathcal{P}$
such that
\begin{eqnarray}
\label{mixed-form}
\begin{cases}
\begin{split}
\int_{\Omega} k^{-1} u \cdot v dx+\int_{\Omega} \nabla p\cdot v dx&=0,  \quad  & \forall v\in \mathcal{V}, \\
\int_{\Omega}\text{div}(u)q dx &=\int_{\Omega} f q dx \quad  & \forall q\in \mathcal{Q},
\end{split}
\end{cases}
\end{eqnarray}
where $\mathcal{U} = H(div, \Omega)$, $\mathcal{P} = H^1(\Omega)$, $\mathcal{V} = (L^2(\Omega))^d$, and
$\mathcal{Q} = L^2(\Omega)$.

There are different ways to construct a discretization of equation \eqref{mixed-form}.
One can try to find the approximations such that
$\{u_h, p_h\} \in H(div, \Omega)\times H^1(\Omega)$ \cite{Thomas_Trujilo-1999}.
Unfortunately these discretizations are computationally expensive and
are only applicable on very restrictive grids. Frequently  dual mixed finite element methods are used
\cite{BreFort91, RobThom91} to construct the discretization $\{u_h, p_h\} \in H(div, \Omega)\times L^2(\Omega)$.
We will consider a class of methods that are related to the primal mixed finite element methods
\cite{RobThom91}, i.e., the discretizations we seek are
$\{u_h, p_h\} \in \mathcal{U}_h \times \mathcal{P}_h, \, \, \mathcal{U}_h \subset \left(L^2(\Omega)\right)^d, \,
\mathcal{P}_h \subset H^1(\Omega)$ with additional conservation enforced on particular volumes.
These methods are close related to the standard cell-centered finite volume method and several
multi-point flux approximation (MPFA) methods that generalize it. We will refer such approximations as mixed finite volume methods.

We assume that two grids are defined:  primary grid $\mathcal{T}_h$  and dual grid $\mathcal{D}_h$.
Usually the primary grid is used to approximate the scaler variable $p$,
and the dual mesh is used to construct the discretization of the velocity.
Different examples of primary and dual grids can be found in \cite{AavBBMann1:98, EdwardsRoger:99, mishev02, mishev03, mc06}.
Consider the discrete problem: \\
\indent
Find $\{u_h, p_h\} \in \mathcal{U}_h \times \mathcal{P}_h$ such that
\begin{eqnarray}
\label{primal_discrete_form}
\begin{cases}
\begin{split}
\int_{\Omega} k^{-1} u_h \cdot v_h dx+\int_{\Omega} \nabla_h p_h\cdot v_h dx&=0,
\quad  \forall v_h\in \mathcal{V}_h, \\
\int_{\Omega}\text{div}_h(u_h)q dx &=\int_{\Omega} f q_h dx. \quad  \forall q_h\in \mathcal{Q}_h,
\end{split}
\end{cases}
\end{eqnarray}
A particular mixed finite volume method will be fully described if we define the grids,
the approximation spaces, and the operators $\nabla_h$, and $\text{div}_h$.
The operator $\text{div}_h$ can be defined in the following way:
\begin{equation} \label{div_h}
\int_{\Omega} \text{div}_h(u_h) q_h \, dx
:= \sum_{V \in \mathcal{T}_h} \sum_{E \in \partial V} \int_{E} u_h \cdot n_E q_h \, ds
\end{equation}
with $\mathcal{Q}_h$ the space of piece-wise constants on the volumes $V$ from the primary grid and $E$ is an
edge in 2-D or face in 3-D in $\partial V$. Here $n_E$ denotes the outward normal vector to $E$.

In the paper, we focus on 2-D case. The 3-D case is a straightforward
extension of 2-D case.

We give an example of mixed finite volume methods as following.

\begin{example} [Mixed FV 1]
We assume that the primary grid $\mathcal{T}_h$ consists of rectangles,
and the dual grid $\mathcal{D}_h$ is also rectangular, but with vertexes the cell centers.
The discrete space  for the scalar variable, $\mathcal{P}_h = \mathcal{Q}_h$, consists
of piecewise constants and is defined on the primary grid,
The approximation space for the vector variable, $\mathcal{U}_h = \mathcal{V}_h$,
is the space of piece-wise constants vector constants with continuous normal components and
is defined on the dual mesh.  Consider one dual cell $D = D_i \cup D_j \cup D_k \cup D_l$
(see Figure \ref{LocalRectangle}). The four functions, $e_{ij}$, $e_{ik}$,
$e_{jl}$ and $e_{kl}$ are defined with the relations
\[
\int_{l_{rt}} e_{pq} \cdot n_{rt} \, dx = \delta _{pq, rt},
\]
where $pq$ and $rt$ can be any element of the set $I_D = \{ ij, \, ik, \, jl, \, jk \}$.
It is easy to see that $\{ e_{pq} \},\, pq \in I_D$  are linearly independent and
therefore form a basis of $\mathcal{U}_{h|D}$.
The degrees of freedom are the integrals of the flux,
i.e, the numbers  $v_{pq} = \int_{l_{pq}} v \cdot n_{pq} \, dx \quad pq \in I_D$.
Then for any $v \in \mathcal{U}_h$, $v_{|D} = \sum_{pq \in I_D} v_{pq} e_{pq}$.
The operator $\nabla_h$ is given by:
\begin{equation} \label{nabla_h}
\int_{\Omega} \nabla_h p_h \cdot v_h \, dx = \sum_{V \in \mathcal{T}_h} \sum_{E \in \partial V}
\int_{E} v_h \cdot n_E \, [p_h]_E \, dx,
\end{equation}
where
\[
[p]_E =
\lim_{
t \rightarrow 0 \atop t >0}
(p(x + t {\bf n}_E) - p(x - t {\bf n}_E))
\]
and the direction of $n_E$ is from left to right or from bottom to top in 2-D.
Note that from the first equation of \eqref{primal_discrete_form} we can express
\begin{equation} \label{u_p relation}
M \vec{u} = \vec{p},
\end{equation}
with $\vec{u}^T = [u_{ij}, u_{ik}, u_{jl}, u_{kl}]^T$,
$\vec{p} = [p_i-p_j, p_i-p_k, p_j-p_l, p_k-p_l]^T$, and $M$ a $4\times 4$ matrix.
We solve for $\vec{u}$ and plug the result in the second equation of
\eqref{primal_discrete_form} to get the final discretization.
\end{example}

For  Example 1, We can rewrite equation \eqref{div_h}
as
\[
 \sum_{V \in \mathcal{T}_h} \sum_{E \in V} \int_{E} u_h \cdot n_E [q_h]_E=-\int_{\Omega} f q_h dx.
\]
which shows that the matrix of the discretization is symmetric. We note that $\mathcal{P}_h$ in  Example 1 is not subspace of $H^1(\Omega)$. Our approximation of $\nabla p$ is nonconforming.

\setlength{\unitlength}{0.008 in}
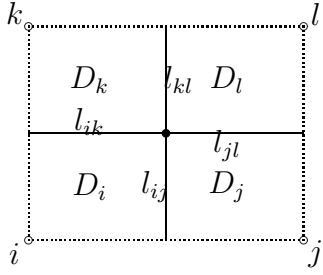
\begin{figure} \centering
\begin{picture}(200,160)
\put(10,10){\circle{5}} \put(190,10){\circle{5}}
\put(10,150){\circle{5}} \put(190,150){\circle{5}}
\put(100,80){\circle*{5}}
\dottedline{3}(10,10)(190,10)(190,150)(10,150)(10,10)
\drawline(100,10)(100,150) \drawline(10,80)(190,80)
\put(1,1){\makebox(0,0){$i$}} \put(1,159){\makebox(0,0){$k$}}
\put(198,1){\makebox(0,0){$j$}} \put(198,158){\makebox(0,0){$l$}}
\put(50,45){\makebox(0,0){$D_i$}}
\put(50,115){\makebox(0,0){$D_k$}}
\put(140,45){\makebox(0,0){$D_j$}}
\put(140,115){\makebox(0,0){$D_l$}}
\put(93,45){\makebox(0,0){$l_{ij}$}}
\put(140,72){\makebox(0,0){$l_{jl}$}}
\put(108,115){\makebox(0,0){$l_{kl}$}}
\put(50,88){\makebox(0,0){$l_{ik}$}}
\end{picture}
\caption{A dual cell $D$}
\label{LocalRectangle}
\end{figure}

The  Example 1 is a straightforward
generalization of the standard cell-centered finite volume method on structured rectangular meshes \cite{mishev03}.
In fact, if the coefficient $k(x)$ in \eqref{model-eq} is a scalar function, then this method
coincides with the standard cell-centered finite volume method.

We need the matrix $M$ in equation \eqref{u_p relation} to have the
appropriate dimensions and to be invertible in order to have a well defined discretization in Example 1.
We state these condition for future reference:
\begin{enumerate} \label{well_defined}
\item $dim(\mathcal{P}_h) + dim(\mathcal{U}_h) = dim(\mathcal{Q}_h) + dim(\mathcal{V}_h)$;
\item Matrix $M$ is invertible.
\end{enumerate}

We will follow Example 1 to derive a new mixed MsFV method in the next section.

\section{A new mixed multiscale finite volume method}
\label{key-mixed-MsFV}


In order to describe  the multiscale finite element method,  we assume that the grids $\mathcal{T}_H$ and $\mathcal{D}_H$
are coarse grids and that there exist an underlying fine grid containing the fine-scale information.
Figure \ref{FVgrid} depicts the rectangle primary coarse grid and dual coarse grid.
The  velocity $u$ is discretized   on the interfaces $\mathcal{E}_H$  of primary grid,
e.g., $E_{ij}$ in  Figure \ref{FVgrid},   and the
pressure $p$ is discretized on cell-centers of the primary mesh / the vertices  of the dual grid,
e.g, vertex $i,j,k,l$ in  Figure \ref{FVgrid}.

Define the discrete space for pressure to be $\mathcal{P}_H = \mathcal{Q}_H$, the space of piecewise constants on $\mathcal{T}_H$,
and for velocity to be $\mathcal{U}_H= \mathcal{V}_H$, a multiscale finite element  space with continuous normal on
$\mathcal{E}_H$ that will be defined
in Subsection \ref{new-Ms-basis}.
Then the mixed multiscale finite volume formulation for  Equation (\ref{mixed-form}) reads:\\
\indent Find $\{u_H, p_H\}\in
\mathcal{U}_H\times \mathcal{P}_H$ such that
\begin{eqnarray}
\label{discrete-form}
\begin{cases}
\begin{split}
\int_{\Omega} k^{-1} u_H\cdot v_H dx +\sum_{E\in \mathcal{E}_H}\int_E v_H\cdot n_E [p_H]_Edx &=0 \ \  \forall v_H\in \mathcal{U}_H,\\
\sum_{V \in \mathcal{T}_H}\int_{\partial K_i} u_H\cdot n q_H dx &=\int_{\Omega} f q_H dx \ \  \forall q_H\in \mathcal{P}_H,
\end{split}
\end{cases}
\end{eqnarray}
where $[p_H]_E$ is the jump of $p_H$ across the interface $E$ and defined in Example 1.
It is clear from equation \eqref{discrete-form} that the operators $\text{div}_H$ and
$\nabla_H$ are defines as follows:
\begin{equation} \label{nabla_H}
\int_{\Omega} \nabla_H p_H \cdot v_H \, dx = \sum_{E\in \mathcal{E}_H}
\int_{E} v_H \cdot n_E \, [p_H]_E \, dx,
\end{equation}
and
\begin{equation} \label{div_H}
\int_{\Omega} \text{div}_h(u_h) q_h \, dx
= \sum_{V \in \mathcal{T}_h} \sum_{E \in \partial V} \int_{E} u_h \cdot n_E q_h \, ds
\end{equation}
Clearly $u_H\in L^2(\Omega)$ and $p_H\in L^2(\Omega)$.
By a similar argument as in Example 1,  the discrete system of (\ref{discrete-form}) is symmetric.

\subsection{A new multiscale velocity basis function}
\label{new-Ms-basis}

In order to complete the derivation of the method from (\ref{discrete-form}), we need to define velocity basis functions.
In this subsection  we design  multiscale velocity basis functions associated with interfaces of $\mathcal{E}_H$.

\begin{figure}[tbp]
\centering
\includegraphics[width=5in, height=3.5in]{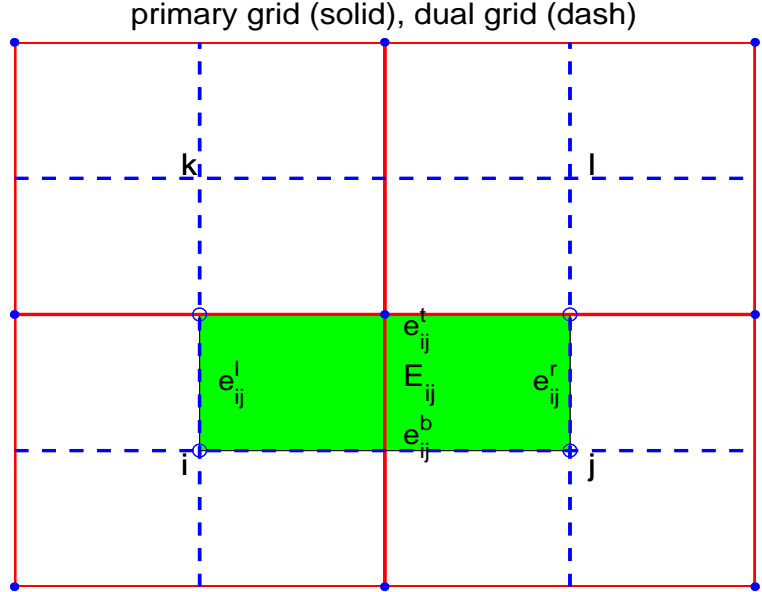}
\caption{Rectangle primary grid and dual grid}
\label{FVgrid}
\end{figure}

\begin{figure}[tbp]
\centering
\includegraphics[width=2 in, height=1.5 in]{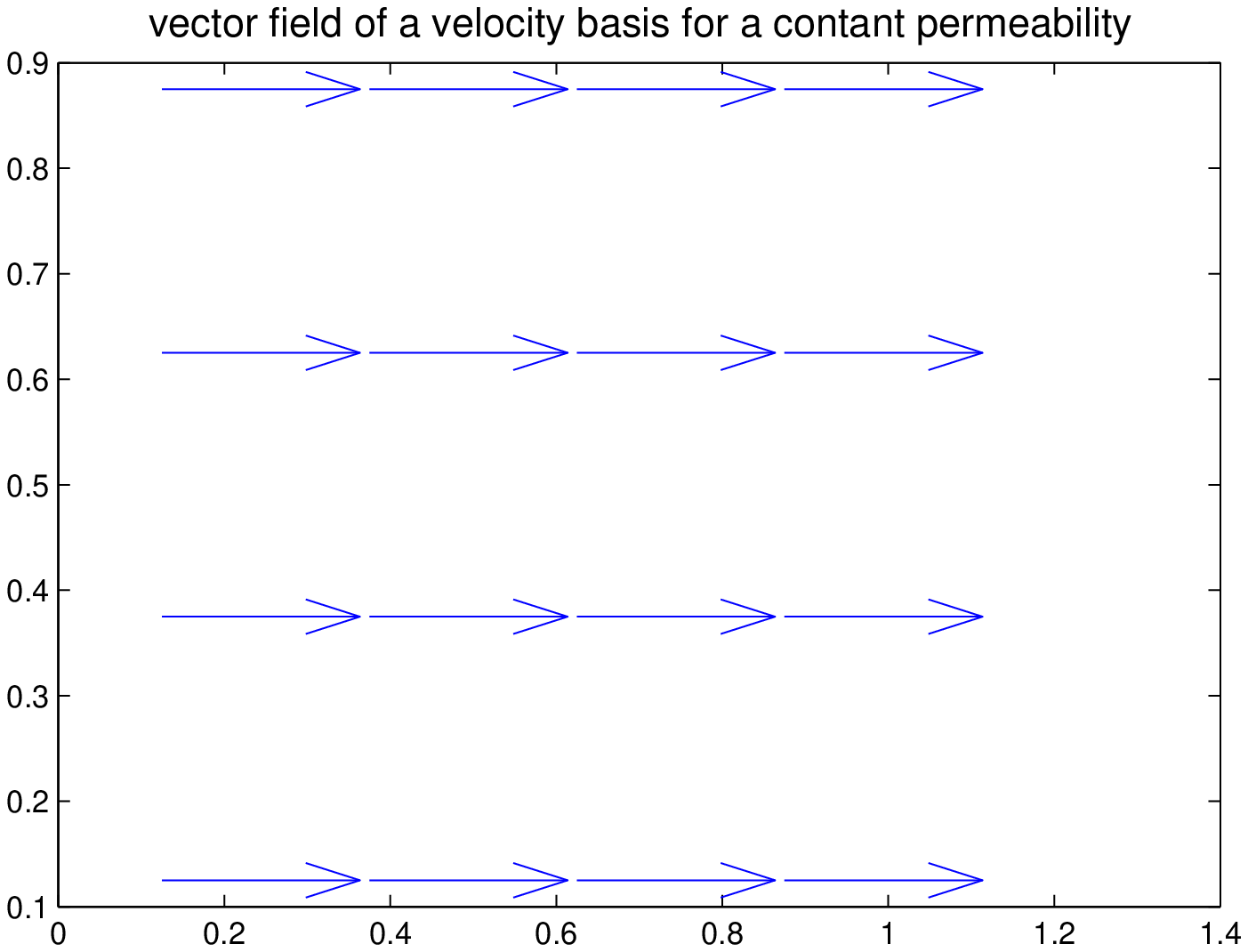}
\includegraphics[width=4 in, height=1.5 in]{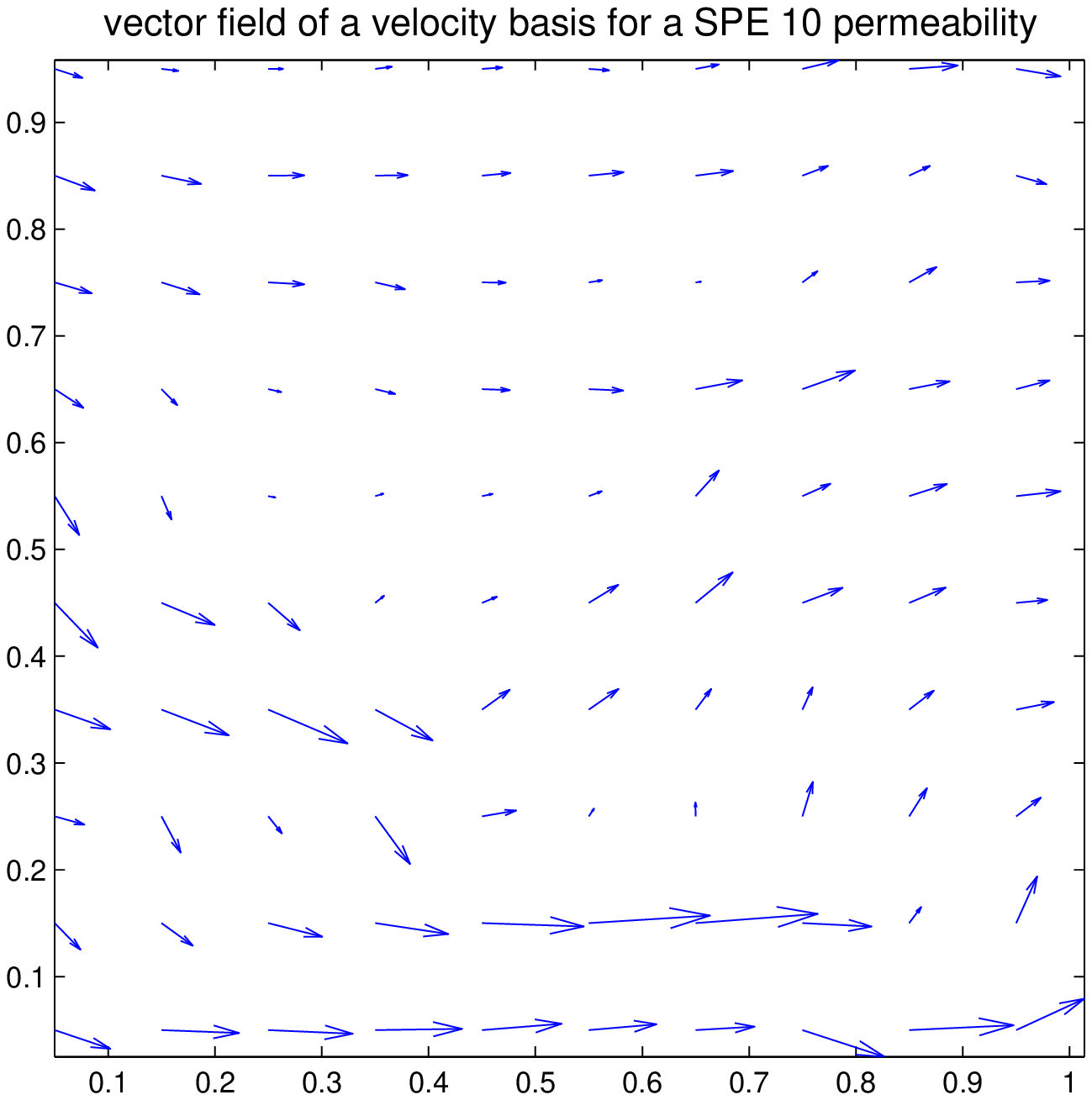}
\caption{Vector field of a velocity basis function for homogeneous  permeability (left)  and
 vector field of a velocity basis function for heterogeneous permeability (right)}
\label{basis-vector}
\end{figure}

  Let $E_{ij}\in \mathcal{E}_H$ be any interface and $\Omega_{ij}$ (green part in Figure \ref{FVgrid}) be an open set bounded by edges  $e_{ij}^l$, $e_{ij}^r$, $e_{ij}^b$ and $e_{ij}^t$.
 We construct a multiscale basis function $\phi_{ij}$ associated to the interface $E_{ij}$ as following:
 \begin{eqnarray}
 \label{basis-equation}
 \begin{cases}
 \begin{split}
 -\text{div}(k\nabla \phi_{ij})&=0  \ \ \text{in} \ \ \Omega_{ij},\\
 -k\nabla \phi_{ij}\cdot n &=\left\{
\begin{array}{ll}
-\frac{v(x)\cdot n}{\int_{e_{ij}^l}v(x)\cdot n dx}    & \text{on} \ \ e_{ij}^l\\
\frac{v(x)\cdot n}{\int_{e_{ij}^l}v(x)\cdot n dx}    & \text{on} \ \ e_{ij}^r\\
0 & \text{on} \ \  e_{ij}^b\cup e_{ij}^t,
\end{array}
\right.
  \end{split}
  \end{cases}
 \end{eqnarray}
 where $v(x)$ is a vector function and has some options depending on the multiscale
 features (e.g., separable scales or non-separable scales).  We will address the options for $v(x)$  later.
  Here $n$ is the unit normal vector pointing out of $\Omega_{ij}$.
  We would like to note that the basis equation (\ref{basis-equation}) is
 defined for a horizontal flux. By switching the no-flow boundary condition and flow boundary condition,
 we can similarly define the basis equation presenting a vertical flux.
  We define the velocity basis function $\psi_{ij}=-k\nabla \phi_{ij}$
 and the finite dimension space for velocity as
 \[
 \mathcal{U}_H=\bigoplus_{E_{ij}\in \mathcal{E}_H} \psi_{ij}.
 \]
 Figure \ref{basis-vector} depicts the vector fields of velocity basis functions (horizontal flux)
 for homogeneous  permeability and heterogeneous  permeability (SPE 10, layer 85), respectively.
 The figure confirms that the multiscale basis defined in (\ref{basis-equation}) reflects the properties
 of the media/permeability. The multiscale basis functions are pre-computed and suited for
parallel computation.

 For many multiscale problems, in particular, the problems with separable scales, we can simply take $v(x)=(1,1)$.  No
 global information is used in this case,  we call this  case the local multiscale
 method.  If the permeability $k$ has strong long range features (e.g.,  highly channelized), we can use
 some global information $v(x)$ related to the features of $u$, e.g.,  the single-phase velocity at time zero,
 to construct the velocity multiscale basis functions \cite {aej08}.
In general, global information used to construct the basis functions in highly heterogeneous
 permeability often   yields much better
approximation  than a constant  boundary condition.  Many numerical studies \cite{aej08, jym09}
show the global information is very helpful to improve accuracy when the permeability
$k$ has distinguished long range features.
We note that one can use time-dependent global fields to construct multiscale basis functions.
This is  helpful for compressible flow simulations.

 If multiple global fields $v_n(x)$ ($n=1,\cdots, N$) are used to build basis functions, then each interface
 corresponds to multiple multiscale basis functions $\phi_{ij}^n$ ($n=1,\cdots, N$), which solve
 \begin{eqnarray}
 \label{basis-multiple}
 \begin{cases}
 \begin{split}
 -\text{div}(k\nabla \phi_{ij}^n)&=0  \ \ \text{in} \ \ \Omega_{ij},\\
 -k\nabla \phi_{ij}^n\cdot n &=\left\{
\begin{array}{ll}
-\frac{v_n(x)\cdot n}{\int_{e_{ij}^l}v_n(x)\cdot n dx}    & \text{on} \ \ e_{ij}^l\\
\frac{v_n(x)\cdot n}{\int_{e_{ij}^l}v_n(x)\cdot n dx}    & \text{on} \ \ e_{ij}^r\\
0 & \text{on} \ \  e_{ij}^b\cup e_{ij}^t.
\end{array}
\right.
  \end{split}
  \end{cases}
 \end{eqnarray}
We note that the multiple global fields $v_n(x)$ can be associated to some representative realizations
in the setting of stochastic two-phase flows \cite{jil10}.

 \begin{rem}
 By the result of Owhadi and Zhang \cite{oz07},   we can use $d$ ($d=dim(\Omega)$) global fields in  (\ref{basis-multiple}).  From the result,  the global field   $v_n(x)=-k\nabla p_n$ ($n=1,\cdots, d$)  are the solutions of the
elliptic equations
\begin{eqnarray}
\label{harmonic-eq}
\begin{cases}
\begin{split}
-\text{div} (k\nabla p_n) &= 0 \ \ \text{in} \ \ \Omega \\
p_n &= x_n \ \ \text{on} \ \ \partial \Omega,
\end{split}
\end{cases}
\end{eqnarray}
where $x=(x_1,\cdots, x_d)$.
\end{rem}

Following the ideas in \cite{jym09}, the global
information can be computed on an intermediate coarse grid  using upscaling techniques. This will
reduce the computation for the global information.

Because the source term of the basis equation (\ref{basis-equation}) is zero,  the flux
conservation implies the following proposition.

\begin{prop}
\label{unit-flux}
Let $E_{ij}\in \mathcal{E}_H$ be any interface and $\psi_{ij}$ be the corresponding
velocity basis function. Then
\[
\int_{E_{ij}}\psi_{ij}\cdot n_{E_{ij}}dx=1.
\]
Moreover, all $\psi_{ij}$ are linear independent, i.e., they form a set of
finite element basis functions.
\end{prop}

\begin{rem}
If $k$ is a constant in $\Omega_{ij}$  and no global information is used, then $\psi_{ij}$ and $\psi_{ik}$ are orthogonal each other,
i.e., the velocity basis function producing horizontal flux is orthogonal to
the velocity basis function producing vertical flux. In fact, we can show
that (see \cite {mishev03})
\[
\psi_{ij}=\left\{
\begin{array}{ll}
{1\over |E_{ij}|}(1,0) & \ \ \text{on} \ \ \Omega_{ij}\\
0  & \ \ \text{else}
\end{array}
\right.
\ \ \text{and} \ \
\psi_{ik}=\left\{
\begin{array}{ll}
{1\over |E_{ik}|}(0,1) & \ \ \text{on} \ \ \Omega_{ik}\\
0  & \ \ \text{else}
\end{array}
\right..
\]
This coincides with the velocity basis in Example 1.
If $k$ is varied in the coarse block,    then $\psi_{ij}$ and $\psi_{ik}$ may not be orthogonal to each other any more.
This can be observed from Figure \ref{basis-vector} (right).
\end{rem}

\subsection{Analysis of the mixed MsFV}
In this subsection we will show that the new mixed MsFV method is well defined.
We have to check the conditions \eqref{well_defined}.
The first one is straightforward since $\mathcal{P}_H = \mathcal{Q}_h$ and $\mathcal{U}_H = \mathcal{V}_H$.
The second condition is verified below. We also give more details how to organize the computations.

The mixed finite volume  formulation
(\ref{discrete-form})  implies the following algebraic linear system
\begin{eqnarray}
\label{PU-system}
\begin{cases}
\begin{split}
\mathbf{A}U + \mathbf{B}P &=0\\
\mathbf{C}U &=F,
\end{split}
\end{cases}
\end{eqnarray}
where $\mathbf{A}$ is a mass matrix.  Here $\mathbf{B}$
has only nonzero entries $1$ and $-1$ and the $\mathbf{B}P$ represents
the jump of $P$. The matrix $\mathbf{C}$  has only nonzero entries $1$ and $-1$, and the sign  depends upon the
normal direction to which  the corresponding flux entries of $U$ associate.
Because  the first equation in (\ref{discrete-form}) can be computed  dual coarse block by dual coarse block
and the support of velocity basis function lies in a dual coarse block, $\mathbf{A}$
can be represented as a block diagonal matrix, i.e.,
 $\mathbf{A}=diag(A_1, A_2,\cdots)$, where  each diagonal
block $A_i$  is the mass matrix associated to a dual coarse block.

A straightforward calculation implies the following proposition.
\begin{prop}
Let $\mathbf{B}$ and $\mathbf{C}$ be defined in (\ref{PU-system}). Then
$\mathbf{B}^T =\mathbf{C}$.
\end{prop}

\begin{figure}[tbp]
\centering
\includegraphics[width=3 in, height=1.5 in]{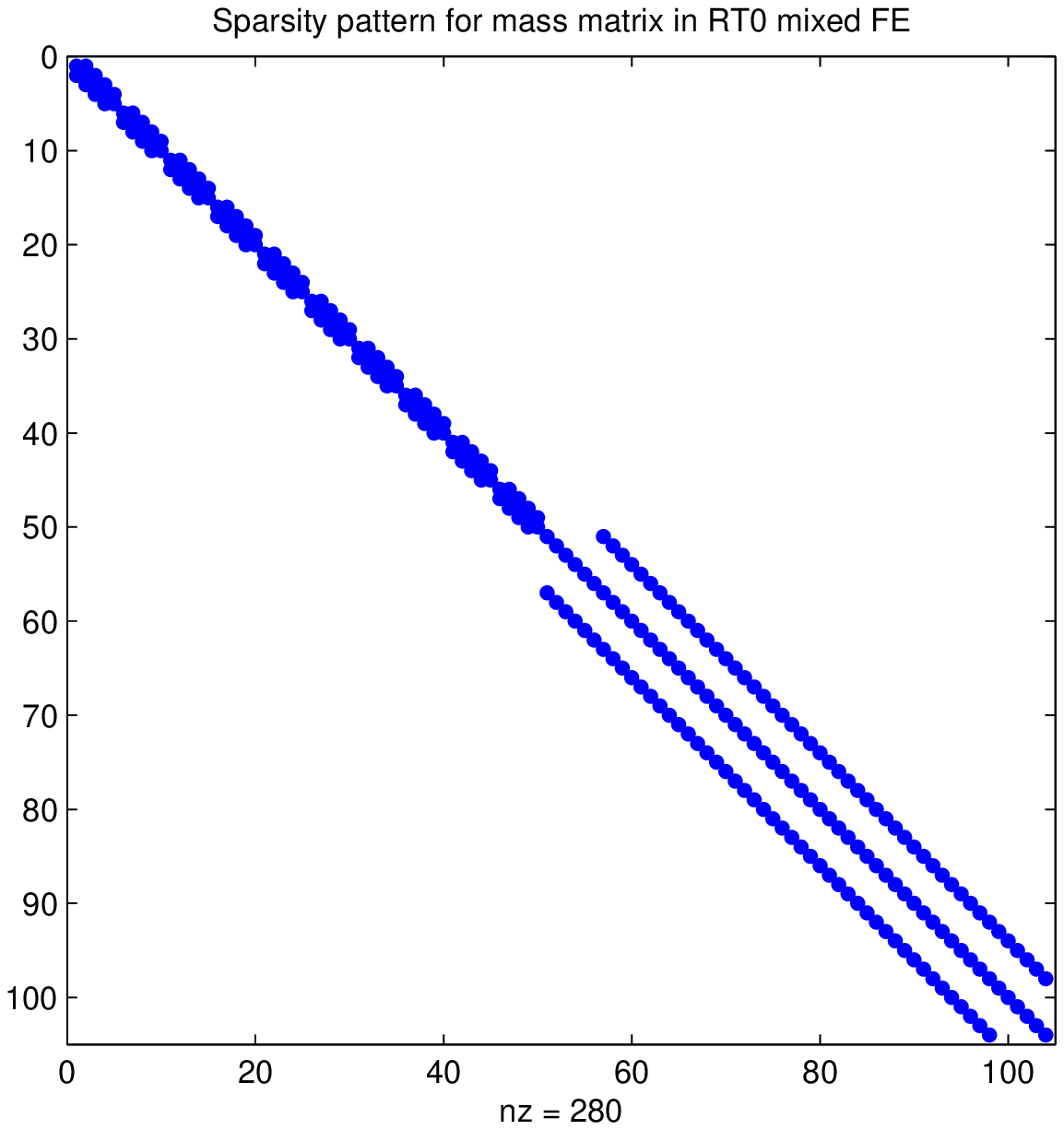}
\includegraphics[width=3 in, height=1.5 in]{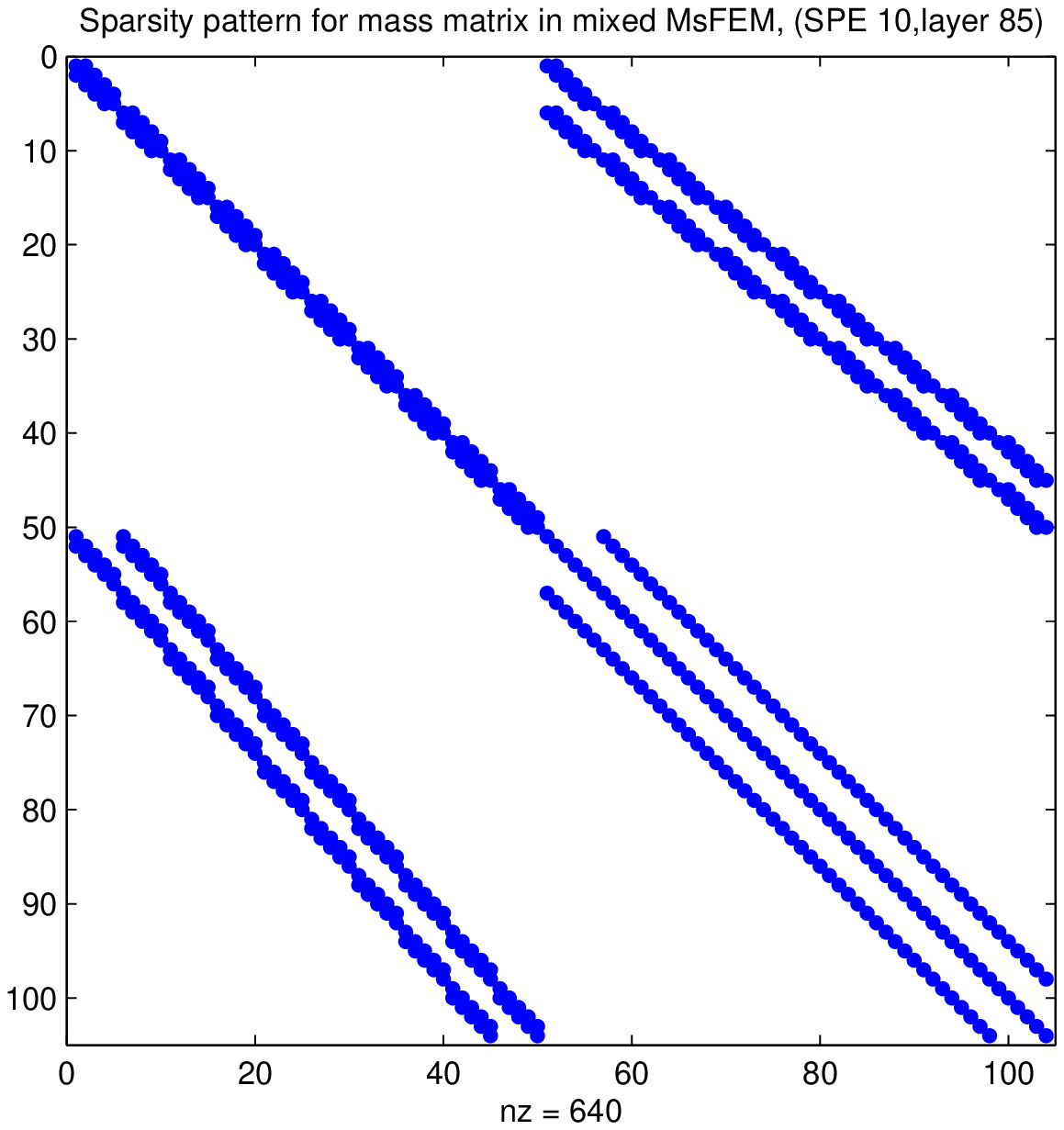}
\caption{Sparsity patterns for  mass matrix in lowest order Raviart-Thomas  mixed FEM (left)  and
 sparsity patterns for  mass matrix in  mixed MsFEM (right)}
\label{mass-matrix}
\end{figure}

Because the mass matrix $\mathbf{A}$ is block diagonal in the mixed MsFV, this allows
to invert each block entry $A_1,A_2,\cdots$ to eliminate the flux $U$ and the computation
(for $A_1^{-1},A_2^{-1},\cdots$) is fast.
It is known that mixed MsFEM also yields  a system such as  (\ref{PU-system}), but
 mass matrix  in mixed MsFEM is not block diagonal and one has to globally compute
the inverse of the mass matrix  to eliminate the flux $U$.   In general, the computation
for the inverse in mixed FEM  is quite time-consuming when the mass matrix  is  large.
This is an advantage of mixed MsFV over mixed MsFEM. Figure \ref{mass-matrix} shows 
the mass matrix for $RT_0$ mixed FEM (left) and mixed MsFEM (right), respectively, where $6\times 10$ grid is used
for both and permeability is heterogeneous (portion  of SPE 10 in layer 85).

We analyze the block diagonal entries of $\mathbf{A}$.  Let $D_{ijkl}$ be the control volume with vertexes  $i,j,k,l$.  In $D_{ijkl}$,  $u_H$
can be represented by
\[
u_H|_{D_{ijkl}}=u_{ij}\psi_{ij}+u_{ik}\psi_{ik}+u_{jl}\psi_{jl}+u_{kl}\psi_{kl}.
\]
Consequently, the first equation in (\ref{discrete-form})
can be reduced to be in $D_{ijkl}$
\begin{equation}
\label{local-eq}
\begin{matrix}
    \begin{bmatrix}
    a_{ij,ij} & a_{ik,ij} & a_{jl,ij} & a_{kl,ij}\\
    a_{ij,ik} & a_{ik,ik}  & a_{jl, ik} & a_{kl, ik} \\
    a_{ij,jl} & a_{ik, jl} & a_{jl, jl}  & a_{kl, jl}\\
    a_{ij, kl}  & a_{ik, kl} & a_{jl, kl}  &  a_{kl, kl}
    \end{bmatrix}
    \begin{bmatrix}
    u_{ij}\\
    u_{ik}\\
    u_{jl}\\
    u_{kl}
    \end{bmatrix}
    +
    \begin{bmatrix}
    p_j-p_i\\
    p_k-p_i\\
    p_l-p_j\\
    p_l-p_k
    \end{bmatrix}
    =
    \begin{bmatrix}
    0\\
    0\\
    0\\
    0
    \end{bmatrix},
    \end{matrix}
\end{equation}
  where $a_{ik,ij}:=\int_{D_{ijkl}}k^{-1} \psi_{ij}\cdot \psi_{ik}dx$ and the other
  entries are defined similarly. We define $A$ to be the most left matrix in  (\ref{local-eq})
and it is a symmetric and positive Gram matrix.  We would like to
 note that $A$ is a representative for the diagonal entries in $\mathbf{A}$ defined in (\ref{PU-system}).
We compute each multiscale basis functions $\psi_{ij}$ in fine scale by standard
mixed FEM (e.g., Raviart-Thomas mixed FEM). Then $\psi_{ij}$  in $D_{ijkl}$ can be represented as
 $\psi_{ij}=\sum_m r_{ij,m} \psi^h_m$,
  where $\psi^h_m$ are standard mixed finite element basis functions in fine scale.
Hence $A$ can be computed in the following way.

  \begin{prop}
  Let $A^h$ be the Gram matrix with entries $(A^h)_{mn}=\int_{D_{ijkl}} k^{-1} \psi^h_{m}\cdot \psi^h_n dx$,
  and each column of $R$ consist of all $r_{ij,m}$.  Then
  \[
  A=R^T A^h R.
  \]
  \end{prop}
  In particular, if permeability  $k$ is homogeneous (constant) in a dual coarse block
  and the mesh is uniformly square and no global information is used,  then a
straightforward calculation implies that the corresponding
 matrix
  \[
  A=diag(2,2,2,2),
  \]
  i.e., $A$ is a diagonal matrix with diagonal entry 2.

Because $A$ is positive (or invertible),  equation (\ref{local-eq}) implies that
\begin{equation}
\label{local-eq2}
\begin{matrix}
     \begin{bmatrix}
    u_{ij}\\
    u_{ik}\\
    u_{jl}\\
    u_{kl}
    \end{bmatrix}
    =
    A^{-1}
    \begin{bmatrix}
    p_i-p_j\\
    p_i-p_k\\
    p_j-p_l\\
    p_k-p_l
    \end{bmatrix}.
    \end{matrix}
\end{equation}
Here $A^{-1}$ is a transmissibility matrix.  We plug the expression (\ref{local-eq2})
into the second equation in (\ref{discrete-form}) to obtain a system  about the pressure,

\begin{equation}
\label{P-eq}
\mathbf{D}P=F.
\end{equation}
Here $\mathbf{D}=\mathbf{C}\mathbf{A}^{-1}\mathbf{B} $, where $\mathbf{B}$
and $\mathbf{C}$ are defined in (\ref{PU-system}). We can show that $\mathbf{D}$
is sparse. In fact, from Equation (\ref{discrete-form}), we know that  only the fluxes of the interfaces
of coarse volume $K_i$ contribute the mass conservation at volume $K_i$,  and that
 each interface flux is determined by its neighbor pressures. Consequently,
there are at most 9 nonzero entries for each row of matrix $\mathbf{D}$ for  $2D$ rectangle cells.
If we do not consider  the restrictions of boundary condition,  $\mathbf{D}$ is symmetric.  This can be shown
by using  the fact $\mathbf{B}^T =\mathbf{C}$.  Our  numerical studies show that
$\mathbf{D}$ is positive.   When permeability is homogeneous, a rigorous mathematics proof
for the positiveness can be found in \cite{mishev03}. Hence  equation (\ref{P-eq}) is solvable.
 Once we obtain the pressure values, then we go back to equation (\ref{local-eq2}) to get the flux
vector $U$.

In particular, if the permeability  $k$ is a constant,  the scheme in (\ref{discrete-form})
coincides with the standard cell-centered difference scheme. In this case, we get the following matrix $\tilde{\mathbf{D}}$ for a $2D$ uniformly square $3\times 3$  grid,
\begin{equation}
\begin{matrix}
\tilde{\mathbf{D}}=
\begin{bmatrix}
2 & -1 & 0 &-1 &0 &0&0 &0 &0  \\
-1&3 &-1 &0 &-1 & 0 & 0 & 0 &0 \\
0&-1 &2 &0 & 0&-1 &0 &0 &0 \\
-1&0&0 &3&-1& 0&-1&0 &0  \\
0 &-1 &0 &-1 &4 &-1 & 0 & -1 &0\\
0 & 0&-1&0 &-1 &3 &0 &0 &-1\\
0 &0 &0 &-1 &0 &0 &2 &-1 &0\\
0 & 0 &0 &0 &-1 &0&-1&3 &-1\\
0& 0& 0&0&0&-1 &0&-1 &2
\end{bmatrix}.
\end{matrix}
\end{equation}
Here we assume Neumann boundary condition without restriction $\int_{\Omega}pdx=0$.
If we restrict $\int_{\Omega}pdx=0$ for a unique solution,
e.g., replace the fifth row in $\tilde{\mathbf{D}}$ by $(1,1,1,1,1,1,1,1,1)$,
the matrix $\tilde{\mathbf{D}}$
becomes an invertible  matrix $\mathbf{D}$.

We summarize the computation as following.\\
{\bf Algorithm 1}
\begin{itemize}
\item For each interface $E_{ij}$, we solve basis equation (\ref{basis-equation}).
\item By equation (\ref{discrete-form}), we formulate an algebraic system (\ref{PU-system}).
\item Eliminate $U$ by local systems  (\ref{local-eq}) and obtain equation (\ref{P-eq}).
\item Solve equation (\ref{P-eq}) to get pressure $P$ and return to equation (\ref{local-eq2}) to get  $U$.
\item By basis equation (\ref{basis-equation}) and $U$, downscale coarse scale velocity to fine scale velocity.
\end{itemize}

\subsection{Reconstruction of the fine-scale velocity field}

Fluxes across the interfaces of  primary coarse grid  can be accurately computed by (\ref{local-eq2}).
In many situations it is often needed to accurately compute the small-scale velocities (or
 fluxes) in regions of interest.  Although we can obtain the velocity in fine grid by straightforwardly
 prolonging  the multiscale basis
function into fine grid,   the  velocity are in general  discontinuous on the
interfaces of the dual coarse grid.  Then large errors can occur in the divergence field and local
mass balance is violated.     We can use
a post-procedure to reconstruct a conservative fine-scale  velocity which is
continuous on the interfaces of both fine grid and coarse grid and fully consistent with
the fluxes across the coarse grid interfaces given by the velocity multiscale basis functions.
We would like to note that a similar post-procedure is used in MsFV \cite{jennylt03}.
We describe the post-procedure as following. \\
{\bf Algorithm 2}
\begin{itemize}
\item Extract the velocities through the interfaces of all coarse
primary cells $K\in \mathcal{T}_H$.  Let $u_H|_{\partial K}$ be
the velocity across $\partial K$, the boundary of $K$. \\
\item On each coarse grid $K$, the following local problem is solved by mixed finite
element methods (or fine volume methods),
\begin{eqnarray}
\label{posterior-local}
\begin{cases}
\begin{split}
-\text{div} (k\nabla \tilde{p})&= f   \ \ \text{in} \ \ K\\
-k\nabla \tilde{p}\cdot n &=u_H|_{\partial K}\cdot n \ \ \text{on} \ \  \partial K.
\end{split}
\end{cases}
\end{eqnarray}
\item Define $\tilde{u}_H= \bigoplus_K \tilde{u}_H^K$, where $\tilde{u}_H^K=-k\nabla \tilde{p}$ and
$\tilde{p}$ solves (\ref{posterior-local}).
\end{itemize}
By the post-procedure, $\tilde{u}_H\in H(div,\Omega)$, that is to say, $\tilde{u}_H$
is continuous along all fine interfaces and coarse interfaces.

\begin{rem}
Because there may exist some errors while computing $u_H$,   $\int_{\partial K} u_H\cdot n dx$
may be  not close to $\int_K f dx$ sufficiently. We can replace the source term $f$ in (\ref{posterior-local})
by ${1\over |K|}\int_{\partial K} u_H\cdot n dx$ to remedy the little disparity for
better accuracy.
\end{rem}
\begin{rem}
It the solution of problem \eqref{posterior-local} is considered to be computationally expensive,
another set of basis functions could be constructed only once and used during the simulation
(see \cite{jennylt03} for a similar approach and  \cite{patent_MMsFV} for details).

\end{rem}

\section {More mixed finite volume methods and their multiscale analogues}
\label{more-FV}

If we define different grids, approximation spaces or operators $\nabla_h$ and $\text{div}_h$,  then we can obtain
different mixed finite volume methods.  One mixed finite volume method has already been described
in Example 1 of Section \ref{key-MFV}, and its multiscale analogue has been developed and analyzed
in Section \ref{key-mixed-MsFV}.   In this section, we briefly  present more examples of
mixed finite volume methods and  their multiscale analogues.

The  following Example 2  is closely related to Example 1.
\begin{example}[Mixed FV 2]
The primary grid $\mathcal{T}_h$ and the dual mesh $\mathcal{D}_h$ are identical to the ones in Example 1.
The discrete space  $\mathcal{P}_h$ for pressure
is the space of bilinear functions. The basis function of $\mathcal{P}_h$ for each cell center is one in the particular
cell center and zero in all neighbors. $\mathcal{Q}_h$ is the space of piecewise constants on the primary grid.
The approximation space for the vector variable, $\mathcal{U}_h = \mathcal{V}_h$
is  the same as in Example 1 and $\nabla_h = \nabla$.
\end{example}

Unstructured grid is often used in practical simulations and
mixed finite volume method can apply to the unstructured grid (see \cite{de06} for extensive discussions).
The following example describe a mixed finite volume method on an unstructured grid.
 \begin{figure}[tbp]
\centering
\includegraphics[width=5in, height=3.5in]{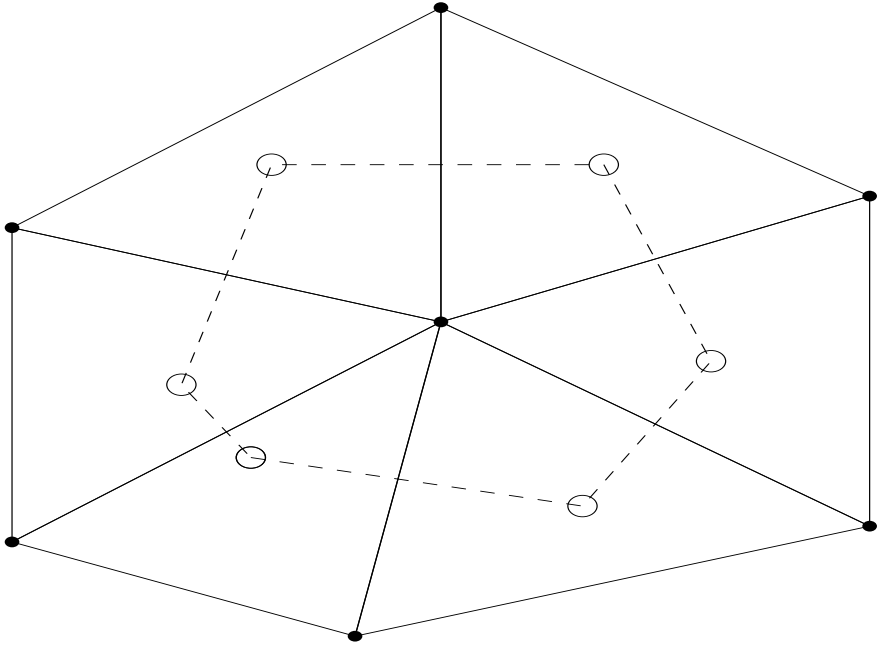}
\caption{Voronoi box/Delaunay triangles}
\label{voronoi}
\end{figure}

\begin{figure}  \centering
\begin{picture}(200,200)
\put(10,10){\circle*{3}}
\put(190,10){\circle*{3}}
\put(120,190){\circle*{3}}
\put(7,0){\makebox(0,0){$(x_i,y_i)$}}
\put(60,0){\makebox(0,0){$e_{ij}$}}
\put(193,0){\makebox(0,0){$(x_j,y_j)$}}
\put(140,0){\makebox(0,0){$e_{ji}$}}
\put(120,200){\makebox(0,0){$(x_k,y_k)$}}
\put(27,66){\makebox(0,0){$e_{ik}$}}
\drawline(10,10)(190,10)(120,190)(10,10)
\put(77,146){\makebox(0,0){$e_{ki}$}}
\put(180,66){\makebox(0,0){$e_{jk}$}}
\put(151,146){\makebox(0,0){$e_{kj}$}}
\put(100,80){\circle{5}}
\put(115,70){\makebox(0,0){$(\bar{x},\bar{y})$}}
\put(65,45){\makebox(0,0){$D_i$}}
\put(140,45){\makebox(0,0){$D_j$}}
\put(113,125){\makebox(0,0){$D_k$}}
\put(107,40){\makebox(0,0){$l_i$}}
\put(79,82){\makebox(0,0){$l_k$}}
\put(125,97){\makebox(0,0){$l_j$}}
\dottedline{3}(100,80)(100,10)
\dottedline{3}(100,80)(65,100)
\dottedline{3}(100,80)(155,100)
\end{picture}
\caption{Triangle $D$}
\label{Delauney}
\end{figure}
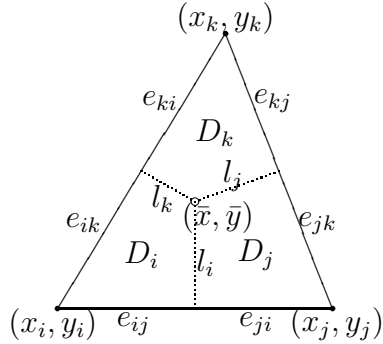
\begin{example}[Mixed FV  3]
We assume that the dual grid $\mathcal{D}_h$ consists of triangles and the cells in the
primary mesh $\mathcal{T}_h$ are control volumes around each vertex in $\mathcal{D}_h$.
Particular example is the Delauney mesh as a dual grid and the corresponding Voronoi grid as primary grid
 (see Figure \ref{voronoi}).
The space $\mathcal{P}_h$ is the space of piece-wise linear functions
on the dual  grid. The space $\mathcal{Q}_h$ is the space of piece-wise constants on the primary grid. The approximation space for the vector variable, $\mathcal{U}_h = \mathcal{V}_h$,
is the space of piece-wise constant vectors with continuous normal components and
it is defined on the dual mesh (see Figure \ref{Delauney}).
The construction of the basis is analogous with the procedure
described in Example 1. The operator $\nabla _h = \nabla$.
Note that the dual grid can be any mesh of triangles, and the primary grid can consists of control
volumes, not necessarily convex polygons, around each vertex. Frequently the control volumes are formed
by connecting the middle of the edges  of each triangle  with the center of mass of the triangle.
The details of the method can be found in \cite{mc06}.

If the dual grid is Delauney mesh and the primary grid is Voronoi mesh, then we can approximate the scalar
variable with piecewise constants, i.e., $\mathcal{P}_h$ is the space of piecewise constants on the primary grid.
Then the operator $\nabla_h$ is defined by  \eqref{nabla_h}.
\end{example}

Here we consider a  mixed finite volume method closely related to the multi-point flux
approximation (MPFA)   discretizations \cite{mishev02}.
\setlength{\unitlength}{.025in}
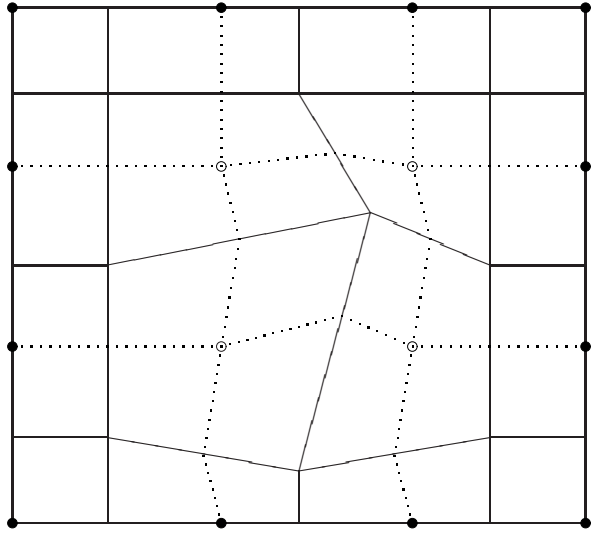
\begin{figure}  \centering
\begin{picture}(140,120)
\drawline(10,10)(130,10)(130,118)(10,118)(10,10)
\drawline(10,28)(30,28)(70,21)(110,28)(130,28)
\drawline(10,64)(30,64)(85,75)(110,64)(130,64)
\drawline(10,100)(30,100)(70,100)(110,100)(130,100)
\drawline(30,10)(30,28)(30,64)(30,100)(30,118)
\drawline(70,10)(70,21)(85,75)(70,100)(70,118)
\drawline(110,10)(110,28)(110,64)(110,100)(110,118)
\put(53.75,47){\circle{2}}
\put(53.75,10){\circle*{2}}
\put(10,47){\circle*{2}}
\put(10,10){\circle*{2}}

\put(93.75,47){\circle{2}}
\put(93.75,10){\circle*{2}}
\put(130,47){\circle*{2}}
\put(130,10){\circle*{2}}

\put(53.75,84.75){\circle{2}}
\put(10,84.75){\circle*{2}}
\put(53.75,118){\circle*{2}}
\put(10,118){\circle*{2}}

\put(93.75,84.75){\circle{2}}
\put(93.75,118){\circle*{2}}
\put(130,84.75){\circle*{2}}
\put(130,118){\circle*{2}}

\dottedline{2}(10,47)(53.75,47)(79,53.4)(93.75,47)(130,47)
\dottedline{2}(10,84.75)(53.75,84.75)(77.5,87.5)(93.75,84.75)(130,84.75)
\dottedline{2}(53.75,10)(50,24)(53.75,47)(57.5,69.5)(53.75,84.75)(53.75,118)
\dottedline{2}(93.75,10)(90,24)(93.75,47)(97.5,69.5)(93.75,84.75)(93.75,118)

\end{picture}
\caption{Qaudrilateral mesh with volumes and covolumes}
\label{quadMesh}
\end{figure}

\setlength{\unitlength}{.010in}
\begin{figure}  \centering
\begin{picture}(300,200)
\drawline(10,10)(20,190)(280,170)(290,15)(10,10)
\dottedline{5}(15,100)(285,92.5)
\dottedline{5}(150,12.5)(150,180)
\put(150,96.25){\circle*{5}}
\put(215,175){\circle*{5}}
\put(282.5,132.25){\circle*{5}}
\dottedline{2}(150,96.25)(215,175)(282.5,132.25)(150,96.25)
\put(17.5,145){\circle*{5}}
\put(12.5,55){\circle*{5}}
\put(85,185){\circle*{5}}
\put(80,11.25){\circle*{5}}
\put(220,13.75){\circle*{5}}
\put(287.5,53.75){\circle*{5}}
\put(80,55) {\makebox(0,0){$1$}}
\put(80,140){\makebox(0,0){$4$}}
\put(210,55){\makebox(0,0){$2$}}
\put(210,135){\makebox(0,0){$3$}}
\end{picture}
\caption{Quadrilateral $V_i$}
\label{Quadri}
\end{figure}
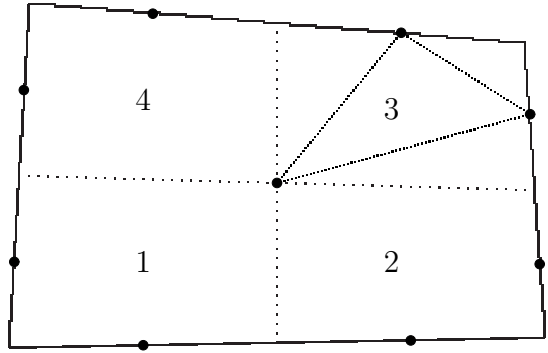

\begin{example}[Mixed FV  4]
The primary mesh is a general quadrilateral mesh (see Fig. \ref{quadMesh}). The dual cells are formed
by connecting the cell-centers with the centers of the edges of the quadrilaterals.
These lines split each primary cell $V_i$ into four quadrilaterals $V_{i,j}, \, j=1,\, 2,\, 3,\, 4$.
The space $\mathcal{P}_h$ is the space of piecewise linear functions on  $V_{i,j}$ with a
common point in the cell center and the other two points on the edges  on the quadrilateral $V_{i,j}$.
(See Fig. \ref{Quadri}).
The space $\mathcal{Q}_h$ is the space of piecewise constant functions on the primary grid.
The space $\mathcal{U}_h$ is the space of piecewise vector constants with continuous normals on the dual cells,
and the space $\mathcal{V}_h$ is the space of piecewise vector constants on the dual cells.
Note that these two spaces have different dimensions.
Each local subspace $\mathcal{U}_{h,i}$ has four degrees of freedom. Each local subspace $\mathcal{V}_{h,i}$
has eight degrees of freedom.  The operator $\nabla_h$ is defined by
\[
\int_{\Omega} \nabla_h p_h \cdot v_h \, dx = \sum_{V_i \in \mathcal{T}_h} \sum_{j=1}^4 \int_{V_{i,j}}
\nabla p_h \cdot v_h \, dx.
\]
Note that $\mathcal{P}_h$ here  is not subspace of $H^1(\Omega)$. The approximation of $\nabla p$ is nonconforming.
One way to make the mixed finite volume  method more robust is to use a conforming approximation.
For example, it is possible to use as $\mathcal{P}_h$ the space of piecewise bilinear functions on  $V_{i,j}$.
One extra basis function is added for each intersection of four adjoint primary cells in two dimensions
\cite{enriched_MPFA}. Extra basis functions are required for three dimensional grids \cite{EdwardsZheng:2008}.

We note that the same procedure works for the grids in Example 3 \cite{VermaAziz:96}.
\end{example}



We used the mixed finite volume methodology to develop a new mixed multiscale method based on
Example 1 in Section \ref{key-mixed-MsFV}. We will sketch how we can follow the same procedure and
derive mixed multiscale finite volume methods based on Example 2 - Example 4.
Some of the methods derived following this framework are known, most are new.

We assume that a consistent mixed finite volume method is defined on the fine mesh, i.e., fine
primary and dual grids are given with the corresponding spaces and the operators.
We need coarse primal and dual grids and we will assume that each coarse cell consists of
fine cell of the same type, i.e., every coarse primary grid cell is a union of adjacent
fine primary grid cells. The same is true for each coarse dual grid cell.
This construction is straightforward for the structured grids.
We suppose that an appropriate coarsening algorithm is used to define the coarse grid for the
unstructured grids. The next step is to define the approximation spaces.
We will require that the coarse discrete spaces provide some approximation of the corresponding functions.
This requirement is easily fulfilled if we follow the same procedure on the coarse structured grids.
For example the space $\mathcal{P}_H$ can consists of constant functions on the primary coarse grid
$\mathcal{T}_H$. Then the space $\mathcal{U}_H$ has to be a multiscale space.
We can consider a  multiscale finite element (e.g., \cite{dmz09, hw97, jennylt03} for detailed description)
space for $\mathcal{P}_H$ and the space $\mathcal{U}_H$
could be the space of piecewise constant vector functions on the dual cells.
The situation is more complicated for unstructured grids. We need approximation of the
gradient of the pressure on unstructured grids and therefore we have to use a multiscale finite element
space $\mathcal{P}_H$. The space $\mathcal{U}_H$ can be either a standard piecewise vector constant space
or a multiscale space.

We provide below more details for several mixed multiscale methods that
can be derived following our methodology and using the mixed finite volume methods on the fine grid in Examples 2 - Example 4.

For better presentation, we call the mixed MsFV proposed in Section \ref{key-mixed-MsFV} to be {\it Mixed MsFV 1}.

\begin{example}[Mixed MsFV 2]
The derivation below is related to Example 2 (Mixed MsFV 2).
The primary grid $\mathcal{T}_H$ and the dual mesh $\mathcal{D}_H$ are identical to the ones in Mixed MsFV 1
(Section \ref{key-mixed-MsFV}).
We chose  $\mathcal{P}_H$ to be the space of multiscale functions. The space $\mathcal{U}_h = \mathcal{V}_h$
is the space of piecewise constant vector functions with continuous normals on the dual coarse grid.
The operator $\nabla_H = \nabla$.

The basis of $\mathcal{P}_H$ on $\mathcal{T}_H$ can be calculated in the same way as in \cite{jennylt03}.
Then we will exactly reproduce the method proposed in \cite{jennylt03}.
We can derive a different method by using a different basis for $\mathcal{P}_H$.
For example, the discrete harmonic basis \cite{XuZikatanov:2004, dmz09} can be constructed in a multilevel way,
that is cheaper, and the computations can be performed using only the fine matrix.

Another method is derived if we select $\mathcal{P}_H$ to be the space of bilinear functions
and $\mathcal{U}_h = \mathcal{V}_h$ to be a multiscale space defined in the same way as in Mixed MsFV1.
The operator $\nabla_H = \nabla$.
\end{example}

The case for unstructured grids is more complicated.
Here we sketch the derivation of a mixed MsFV  for unstructured fine triangular
grids. 


\begin{example}[Mixed MsFV 3]
We assume that the primary and dual grids $\mathcal{T}_H$, $\mathcal{D}_H$
are constructed from fine cells described in Example 3 (MFV 3) using an appropriate coarsening algorithm \cite{XuZikatanov:2004}.
The first method we propose uses a multiscale space $\mathcal{P}_H$ of discrete harmonic functions
discussed in previous example. The space  $\mathcal{U}_H = \mathcal{V}_H$ is the space of
constant vector functions with continuous normals.
The operator $\nabla_H=\nabla$.
Note that this method does not use global information.

We recommend using the multiscale space $\mathcal{U}_H = \mathcal{V}_H$ when it is beneficial to use some global information about the problem and
this information can be transferred using the velocity. The basis functions for $\mathcal{U}_H$ are
constructed in a similar way as on the rectangular grid.
We modify the definition of the velocity basis function
(\ref{basis-equation}) as follows. Consider the triangle $D$ on Figure  \ref{Delauney}.
The basis function corresponding to $l_i$ is the solution of the following boundary value problem:
 \begin{eqnarray}
 \label{basis-equation-unstr}
 \begin{cases}
 \begin{split}
 -div(k\nabla \phi_{ij})&=0  \ \ \text{in} \ \ D_i \cup D_j,\\
 -k\nabla \phi_{ij}\cdot n &=\left\{
\begin{array}{ll}
-\frac{v(x)\cdot n}{\int_{e_{ik}}v(x)\cdot n dx}    & \text{on} \ \ e_{ik}\\
 \frac{v(x)\cdot n}{\int_{e_{jk}}v(x)\cdot n dx}    & \text{on} \ \ e_{jk}\\
0                                                     & \text{on} \ \ l_{k} \cup l_{j}\\
0                                                     & \text{on} \ \ e_{ij} \cup e_{ji}.
\end{array}
\right.
  \end{split}
  \end{cases}
 \end{eqnarray}
Again $\nabla_H=\nabla$.
\end{example}

\begin{example}[Mixed  MsFV 4]
There are several ways to construct multiscale MPFA methods.
One approach is presented in \cite{jennylt03}.
This requires a construction of a multisclale space $\mathcal{P}_H$ and the standard spaces for
$\mathcal{U}_H$ and $\mathcal{V}_H$. Another method with the ability to utilize the available global information
can be constructed by using a multiscale space $\mathcal{P}_H$, not necessarily the same as in \cite{jennylt03}, and
a multiscale spaces $\mathcal{U}_H$ and $\mathcal{V}_H$ defined in similar way as  \eqref{basis-equation-unstr}.
There exists a third way for rectangular or quadrilateral grids, such that
the edges  of the coarse cells are straight lines.
We can use the standard finite element space $\mathcal{P}_H$ and a multiscale
space $\mathcal{U}_H$. The framework also can be applied to the MPFA discretizations proposed in
\cite{enriched_MPFA, EdwardsZheng:2008, VermaAziz:96} and the corresponding mixed MsFV methods derived.
\end{example}

We note that it  is difficult to apply standard multiscale basis functions in production code because of the
geometric grid information necessary to impose the appropriate boundary conditions. It is also  difficult to apply them  for unstructured grids.
If the coarse grid is very distorted, it may happen that a point where the four quadrilaterals meet is not in the
support of the neighboring functions. This will decrease the accuracy.
If discrete harmonic and multilevel basis is used, then computation becomes cheaper and
 can be applied on the matrix level (see \cite{dmz09}). The extensive study for these issues
  is still  under investigation.

\section{Numerical results}

In this section, we apply the  mixed MsFV propped in Section \ref{key-mixed-MsFV} to
simulate  incompressible two-phase flows in  porous media.
We will consider the porous media with non-separable scales and separable scales.
In the first example, the permeability field is   from
SPE Comparative Solution Project \cite{cb01}
(also known as SPE 10) and its scales are non-separable, and it has channelized structure.
In the second example, we consider the flows in two-point correlation permeability.
The permeability is described with a two-point correlation function and
has non-separable scales and distinct spatial variation.  In the third example,
the permeability is described by a periodic function with a small period
and its scales are apparently separable.  We apply the local mixed MsFV
and global mixed MsFV to the flows in the three types of permeability fields.
We will find that the mixed MsFV can provide accurate approximation
on coarse grid and that using global information is able to greatly improve
accuracy for the cases of non-separable scales.

In our numerical simulations, we will perform two-phase flow and transport
simulations. The equations are given (in the absence of gravity and
capillary effects) by
flow equations
\begin{equation}
\text{div}(\lambda(S) { k} \nabla p)=f, \label{fspeqn}
\end{equation}
where the total mobility $\lambda(S)$ is given by
$\lambda(S)=\lambda_w(S)+\lambda_o(S)$ and $f$ is a source term.
Here, $\lambda_w(S)=k_{rw}(S)/\mu_w$ and  $\lambda_o(S)=k_{ro}(S)/\mu_o$
where $\mu_o$ and  $\mu_w$ are viscosities of oil and water phases, correspondingly,
and $k_{rw}(S)$ and $k_{ro}(S)$ are relative permeability of oil and water phases, correspondingly.
The saturation is governed by
\begin{equation}
{\partial S\over\partial t}+ \text{div}( { F} )=0, \label{fsseqn}
\end{equation}
where ${ F}={ v}f_w(S)$, with $f_w(S)$,
the fractional flow of water, given by $f_w=\lambda_w /
(\lambda_w+\lambda_o)$, and the total velocity ${ v}$ by:
\begin{equation}
{ v} ={ v}_w + { v}_o=-\lambda(S) { k } \nabla p. \label{tveqn}
\end{equation}
In our simulations, we take $k_{rw}(S)=S^2$ and $k_{ro}(S)=(1-S)^2$.
In the presence of capillary effects, an additional diffusion term is present
in (\ref{fsseqn}) and an efficient  treatment of capillarity is proposed in \cite{ks10}.
We note that the porosity is 1  in the saturation equation  (\ref{fsseqn}).

We solve the two-phase flow system (\ref{fspeqn}) and (\ref{fsseqn}) by the classical IMPES (implicit pressure and explicit saturation). The saturation equation (\ref{fsseqn}) is discretized  in fine grid by upwind finite volume method. The temporal discretisation is an implicit scheme, which is unconditionally stable but produce a nonlinear system (Newton-Raphson
iteration solves the nonlinear system). For completeness,  we describe the upwind finite volume method
for equation (\ref{fsseqn})  in Appendix \ref{app1}.

We compare the saturation fields and
water-cut data as a function of pore volume injected (PVI).
The water-cut  is defined as the fraction of water in
the produced fluid and is given by $q_w/q_t$, where $q_t=q_o+q_w$,
with $q_o$ and $q_w$ being the flow rates of oil
and water at the production edge of the model.
In particular, $q_w=\int_{\partial \Omega^{out}} f(S)
{ v}\cdot { n} ds$,
$q_t=\int_{\partial \Omega^{out}} {  v}\cdot { n} ds$,
where ${\partial \Omega^{out}}$ is the outer flow boundary.
Pore volume injected, defined as $PVI={1 \over V_p} \int_0^t
q_t(\tau) d\tau$, with $V_p$ being the total pore volume of
the system, provides the dimensionless time for the displacement.
 We consider a traditional quarter
five-spot
problem, where the water is injected at
left bottom  corner
and oil is produced at  the right top corner of the rectangular domain.
In all numerical simulations, multiscale basis functions are constructed
once at the beginning of the computations. In the discussions,
we refer to the grid where multiscale basis functions are constructed
as a coarse grid. We use the global single-phase information (where $\lambda(S)=1$)
to construct mixed MsFV basis functions.  The global information is computed on fine grid
at time zero and will not change throughout the simulation.

In the simulations, we solve the pressure equation on the coarse grid by {\it Algorithm 1}
and use the post-procedure described in {\it Algorithm 2}
to  re-construct the fine-scale velocity field which is used to solve the saturation equation.

We solve the two-phase pressure equation (\ref{fspeqn}) by the mixed MsFV (mixed FEM for reference solution).
For the numerical simulations, we use 10 time steps for pressure equation, and for each pressure time
step, we use 10 time steps to solve saturation equation. Hence, the time step for
pressure is 0.1 PVI and the time step for saturation is 0.01 PVI.

To assess the performance of the saturations and water-cuts obtained using the  mixed MsFV,  we compute the time-dependent pressure equation on fine grid
 by using lowest order Raviart-Thomas   mixed finite element method, and this produces
a reference velocity to solve a reference saturation solution $S_{ref}$.  By the reference saturation and the reference velocity, we
get the reference water-cut $W_{ref}$.  We measure the relative saturation error in $L^1$-norm
and the relative water-cut error in $L^2$-norm,
\[
\|S_{MsFV}-S_{ref}\|_{L^1}/\|S_{ref}\|_{L^1},  \quad  \|W_{MsFV}-W_{ref}\|_{L^2}/\|W_{ref}\|_{L^2}.
\]
where $S_{MsFV}$ and $W_{MsFV}$ denote the saturation and water-cut  by the mixed MsFV, respectively.

\subsection{Flow in SPE 10 permeability}

\begin{figure}[tbp]
\centering
\includegraphics[width=5in, height=2in]{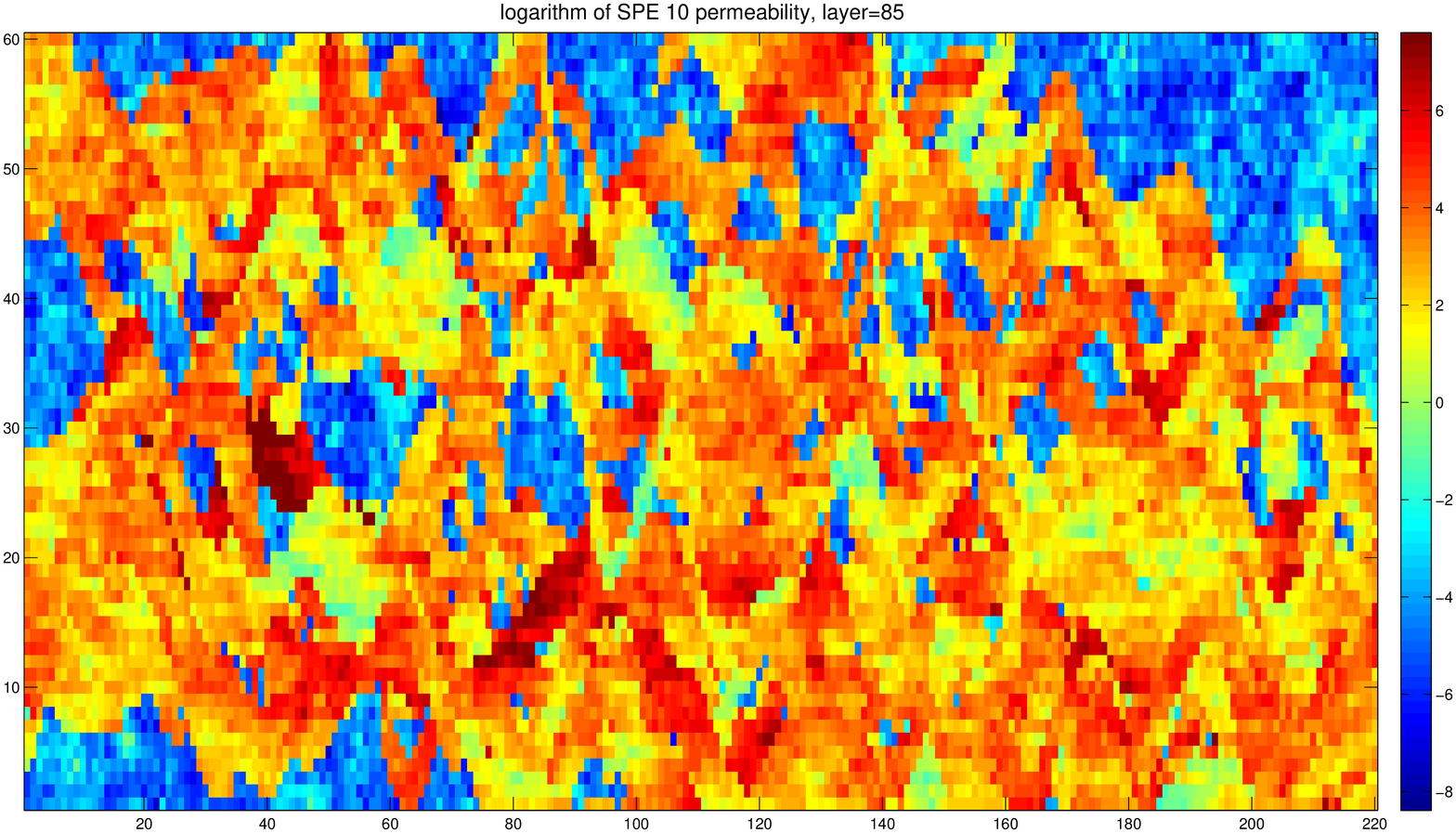}
\caption{ Logarithm of SPE 10 permeability, layer 85.}
\label{exam1-perm}
\end{figure}

For the first  numerical example, we choose the SPE 10 permeability (layer 85), which is highly channelized
and defined on a $60\times 220$ find grid.
The permeability map is depicted in Figure \ref{exam1-perm}. We take $3\times 5$ coarse grid  for
 both the global mixed MsFV and the local mixed MsFV.  We take tests for two different viscosity  ratios of water and oil.

\begin{figure}[tbp]
\centering
\includegraphics[width=6in, height=4in]{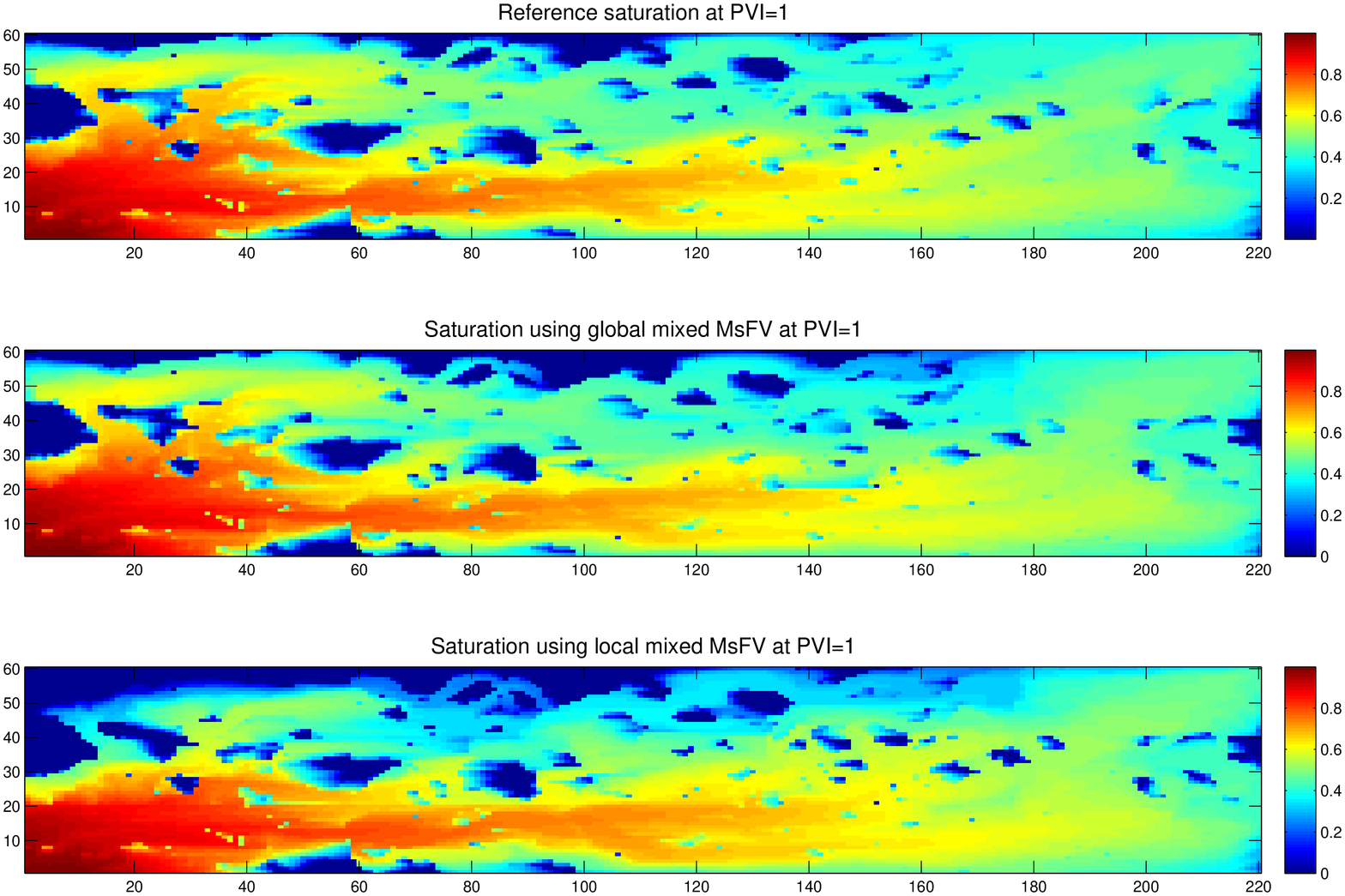}
\caption{Saturation profiles, $\mu_w/\mu_0=1/10$.  Top: Reference saturation at PVI=1.
Middle: Saturation at PVI=1 by the global mixed MsFV. Bottom: Saturation at PVI=1 by the local mixed MsFV.}
\label{exam1-saturation}
\end{figure}

\begin{figure}[tbp]
\centering
\includegraphics[width=4in, height=2in]{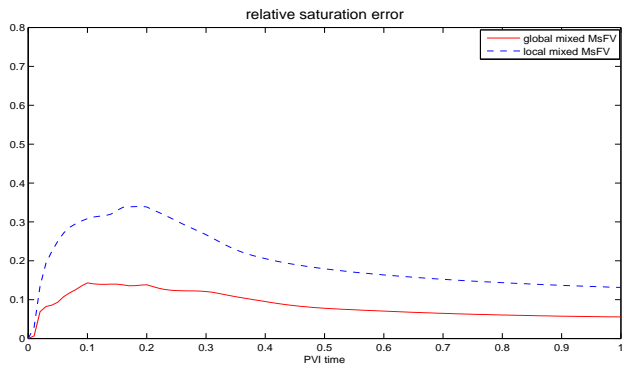}
\caption{saturation error via time, $\mu_w/\mu_0=1/10$.}
\label{exam1-saturation-error}
\end{figure}

\begin{figure}[tbp]
\centering
\includegraphics[width=4in, height=2in]{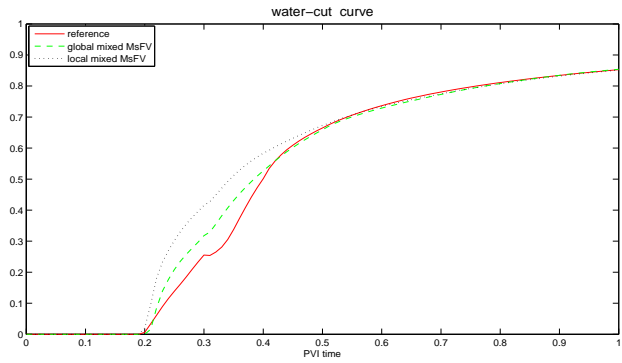}
\caption{water-cut curves, $\mu_w/\mu_0=1/10$.}
\label{exam1-WC}
\end{figure}

 For the first test of the example, we consider the case that viscosity  ratio  $\mu_w/\mu_o$ of water and oil is less than 1.
 Here we take $\mu_w/\mu_o=1/10$ for the simulation.
Figure \ref{exam1-saturation} shows
the reference (fine-scale) saturation profile, the saturation profile
 using the global
mixed MsFV and the saturation profile using the local mixed MsFV, respectively, at PV1=1.
   Figure \ref{exam1-saturation-error} shows the saturation error
via different PVI times.  From the figure, we find that the mixed MsFV using
global information render more accurate saturation solutions throughout the whole simulation than the
saturation from the local mixed MsFV.
The water-cut curves are shown in Figure \ref{exam1-WC}
for reference, global mixed MsFV and local mixed MsFV, respectively.
The water break-through time
is almost the same for the three methods, however,
the water-cut curve by using global mixed MsFV
is closer to the reference water-cut curve at early time than the local mixed MsFV.
The average errors of saturation and water-cut
are shown in Table \ref{exam1-tab}.  
We observe that the error of the solution of the global mixed MsFV is at least two times smaller than the
error of the solution of the local MsFV.

\begin{figure}[tbp]
\centering
\includegraphics[width=6in, height=4in]{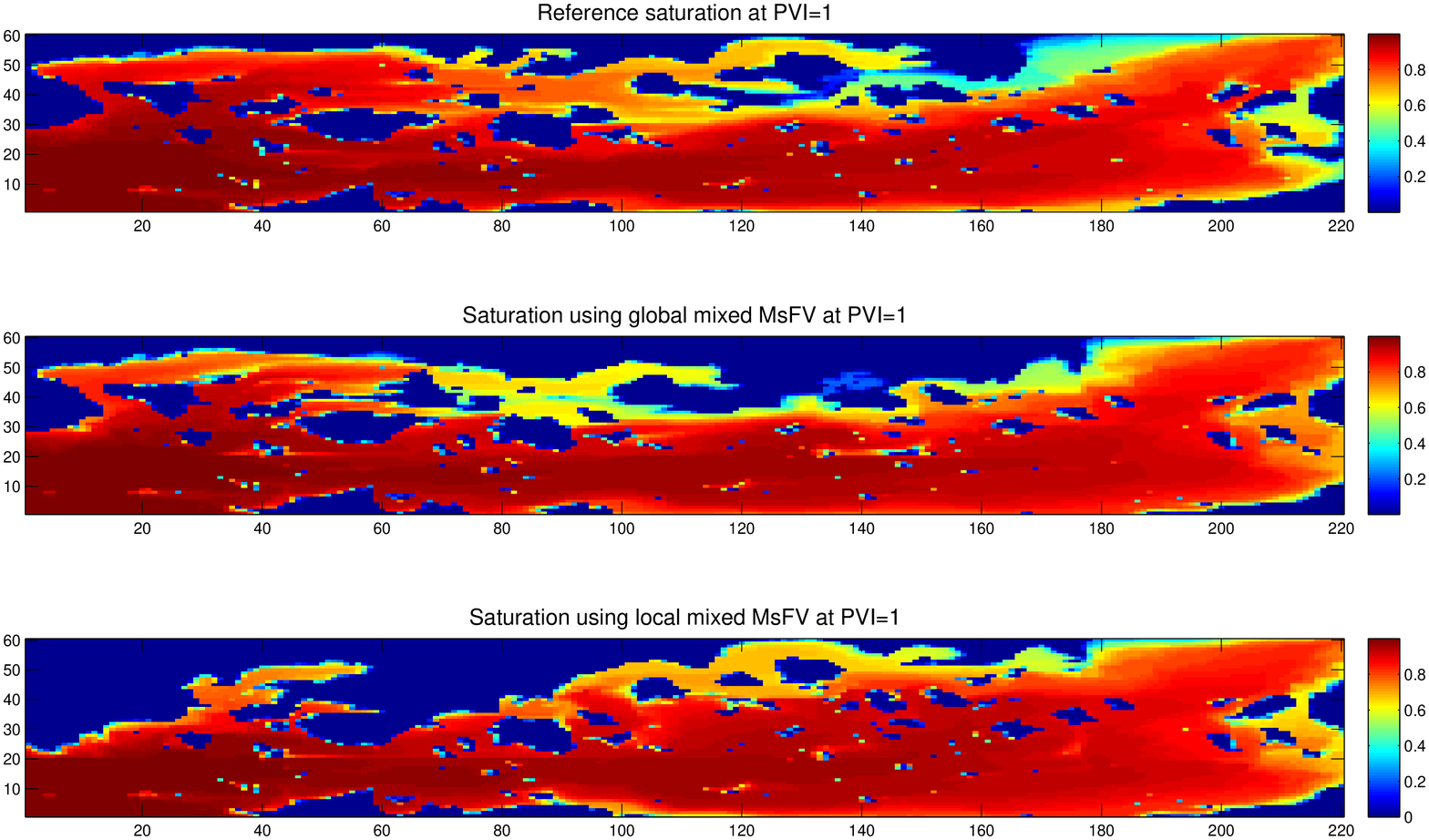}
\caption{Saturation profiles, $\mu_w/\mu_0=3$.  Top: Reference saturation at PVI=1 .
Middle: Saturation at PVI=1 by the global mixed MsFV. Bottom: Saturation at PVI=1 by the local mixed MsFV.}
\label{exam1-2ndtest-saturation}
\end{figure}

\begin{figure}[tbp]
\centering
\includegraphics[width=6in, height=2in]{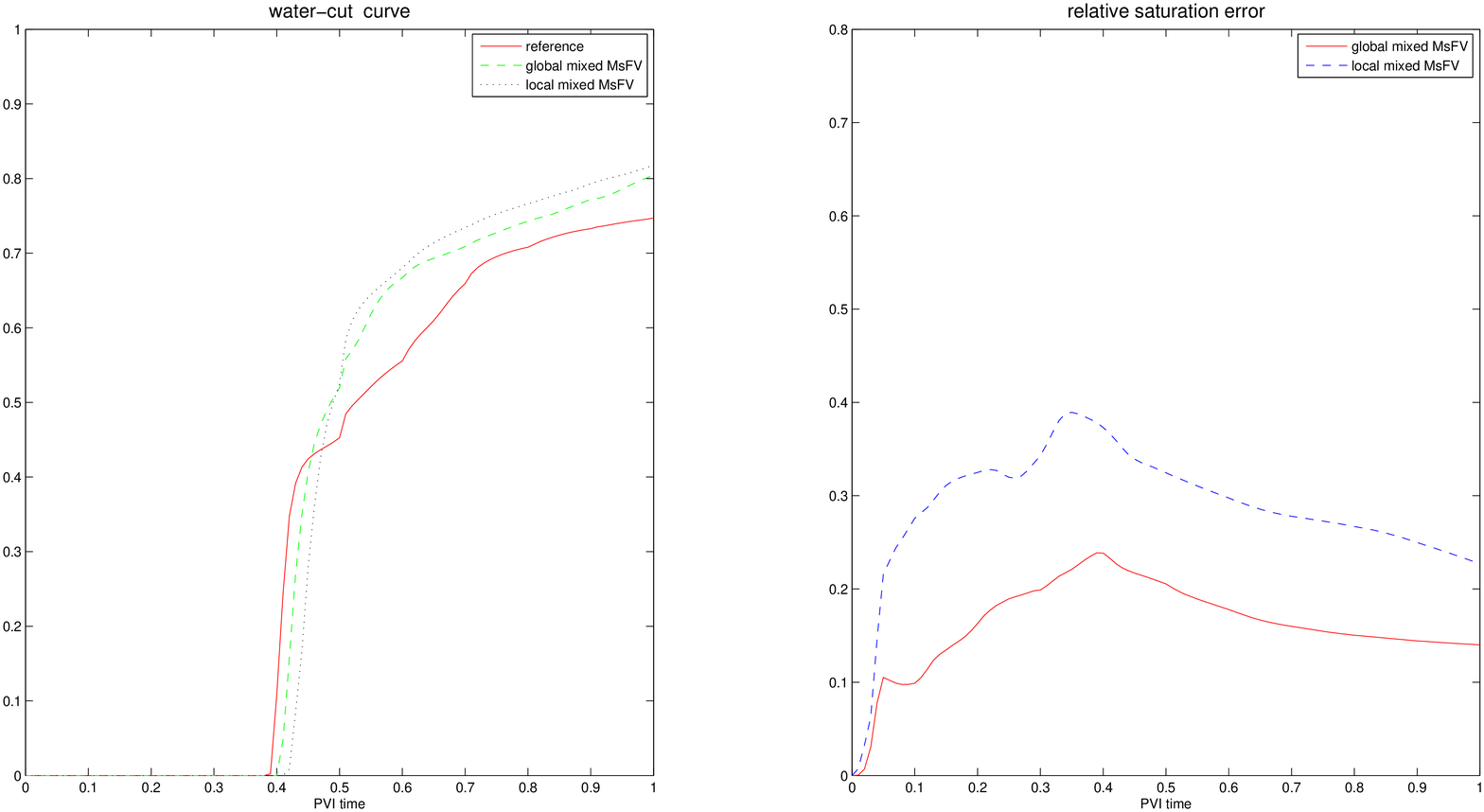}
\caption{water-cut curves and saturation error via time,  $\mu_w/\mu_0=3$.}
\label{exam1-2ndtest-ws}
\end{figure}

For the second test of the example, we consider the case when $\mu_w/\mu_o=3$.
Figure \ref{exam1-2ndtest-saturation} depicts
the reference  saturation profile and  the saturation profiles uisng both  the global
mixed MsFV and  the local mixed MsFV at PV1=1.
Figure \ref{exam1-2ndtest-ws} illustrates  the saturation error
via different PVI times  and the water-cut curves as well.
From Figure \ref{exam1-2ndtest-saturation} and Figure \ref{exam1-2ndtest-ws}, we find that global mixed MsFV performs better than local mixed MsFV in the point of view
of accuracy and water-breakthrough time. Table \ref{exam1-tab} shows  the average errors of saturation and water-cut for the test.

\begin{table}[hbtp]
\centering \caption{relative errors for saturation and water-cut,  $\mu_w/\mu_0=1/10$ }
\begin{tabular}{|c|c|c|c|c|}
\hline
mixed MsFV  &  Water-Cut Error   &  Saturation Error \\

\hline
local mixed MsFV   & 0.1140 & 0.2024 \\
\hline
global mixed MsFV  & 0.0510  & 0.0870\\
\hline
\end{tabular}
\label{exam1-tab}
\end{table}

\begin{table}[hbtp]
\centering \caption{relative errors for saturation and water-cut,  $\mu_w/\mu_0=3$ }
\begin{tabular}{|c|c|c|c|c|}
\hline
mixed MsFV  &  Water-Cut Error   &  Saturation Error \\

\hline
local mixed MsFV   & 0.1748 & 0.2892 \\
\hline
global mixed MsFV  & 0.1062  & 0.1535\\
\hline
\end{tabular}
\label{exam1-2ndtest-tab}
\end{table}

\subsection {Flow in two-point correlation permeability}
\begin{figure}[tbp]
\centering
\includegraphics[width=5in, height=2in]{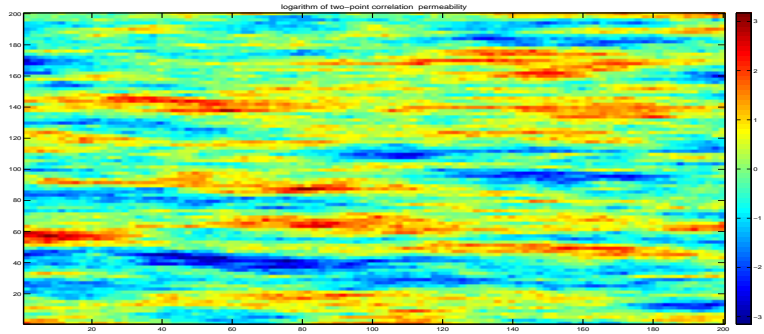}
\caption{ Logarithm of two-point correlation  permeability}
\label{exam2-perm}
\end{figure}

\begin{figure}[tbp]
\centering
\includegraphics[width=6in, height=4in]{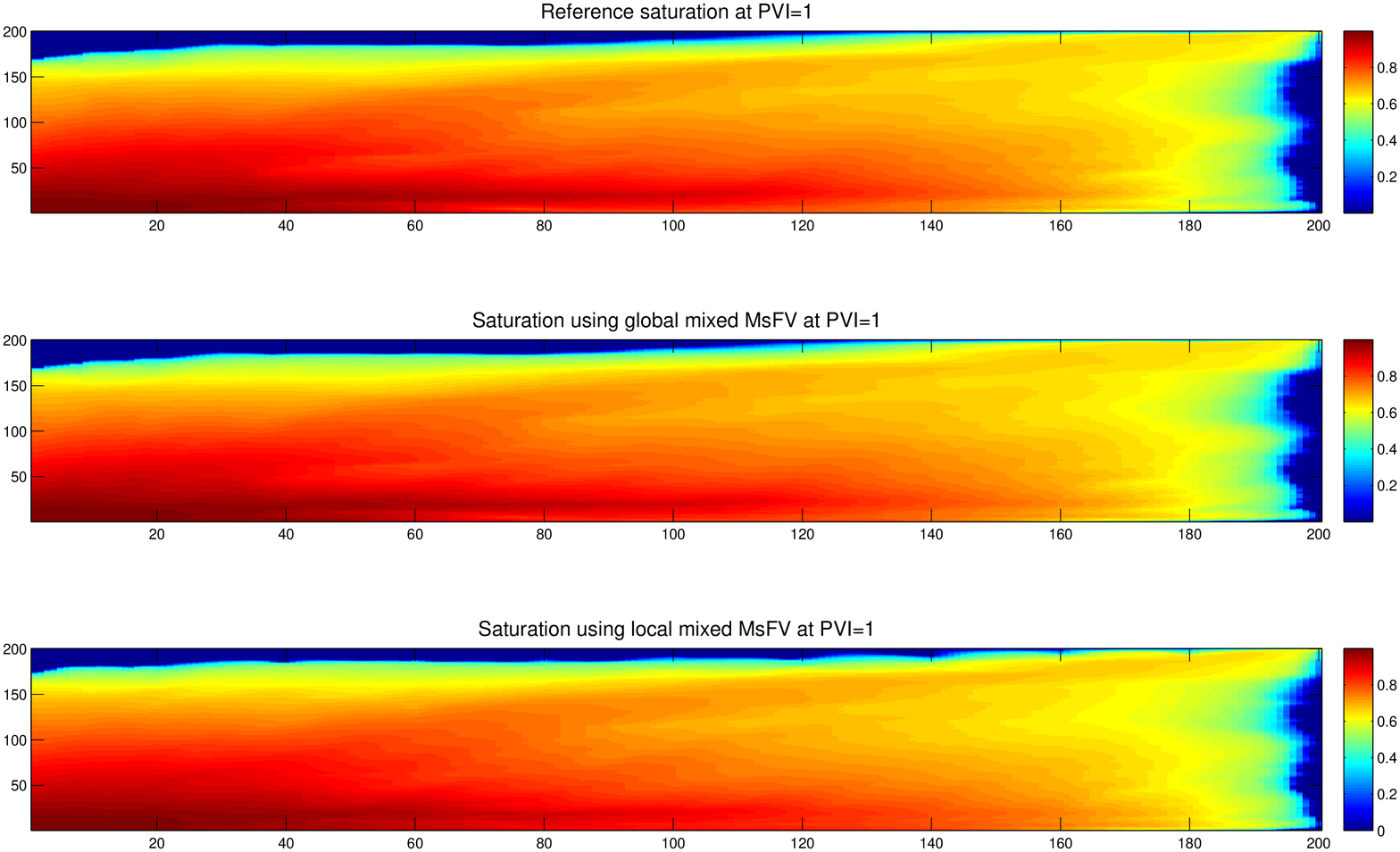}
\caption{Top: Reference saturation at PVI=1 for the second example.
Middle: Saturation at PVI=1 by the global mixed MsFV. Bottom: Saturation at PVI=1 by the local mixed MsFV.}
\label{exam2-saturation}
\end{figure}

\begin{figure}[tbp]
\centering
\includegraphics[width=4in, height=2in]{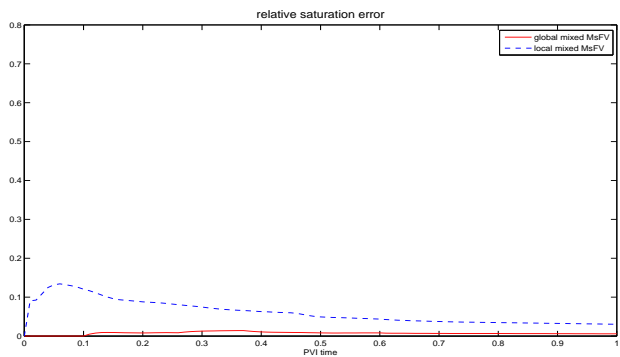}
\caption{saturation error via time for the second example}
\label{exam2-sat-error}
\end{figure}

\begin{figure}[tbp]
\centering
\includegraphics[width=4in, height=2in]{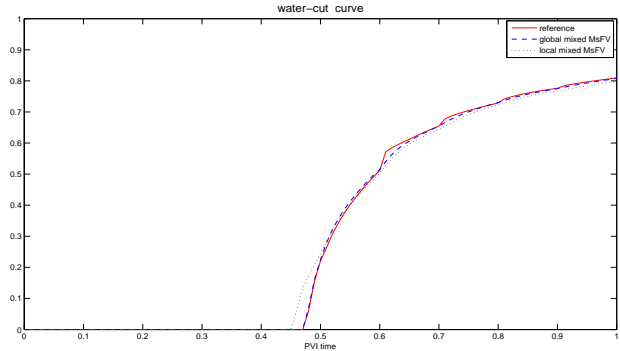}
\caption{water-cut curves for the second example}
\label{exam2-wc}
\end{figure}

In the second example,
 we choose  a realization of the
permeability field generated using  a two-point correlation
function with correlation lengths in $x_1$-direction $L_1=0.4$
and in $x_2$-direction $L_2=0.05$.   Exponential variogram is selected
(see e.g., \cite{dj98}) to generate the permeability. The permeability is defined on
 $200\times 200$ fine-grid and depicted in Figure \ref{exam2-perm}.
The viscosity ratio is $\mu_w/\mu_o=1/3$ and the mixed MsFVs are implemented
on $10\times 10$ coarse grid.
Figure \ref{exam2-saturation} depicts the
reference  saturation, the saturation
field using the global
mixed MsFV and the saturation field using the local mixed MsFV, respectively.
   Figure \ref{exam2-sat-error} demonstrate  the relative  saturation error
via different PVI times.  From the figure, we observe  that mixed MsFV
can provide a good approximation for the flow, and that using global information
improves the accuracy.
The water-cut curves are depicted  in Figure \ref{exam2-wc}
for reference, global mixed MsFV and local mixed MsFV, respectively.
Figure \ref{exam2-wc} shows that the water-cut curve  in global mixed MsFV
is more close to the reference water-cut curve   than the mixed MsFV without
using global information, and that the water break-through time in global mixed MsFV is almost the same as the reference
water break-through time.    The average errors of saturation and water-cut
are shown in Table \ref{exam2-tab}.  From Table \ref{exam2-tab},  we can observe:
 the saturation in the global mixed MsFV is almost 9 times better   than
the saturation in the local mixed MsFV, and the water-cut in global mixed MsFV is around 3 times better
than the water-cut in the local mixed MsFV.

\begin{table}[hbtp]
\centering \caption{relative errors for saturation and water-cut for the second example}
\begin{tabular}{|c|c|c|c|c|}
\hline
mixed MsFV  &  Water-Cut Error   &  Saturation Error \\

\hline
local mixed MsFV   & 0.0464 & 0.0610 \\
\hline
global mixed MsFV  &     0.0127  &     0.0075\\
\hline
\end{tabular}
\label{exam2-tab}
\end{table}

\subsection {Flow in periodic permeability}

\begin{figure}[tbp]
\centering
\includegraphics[width=5in, height=2in]{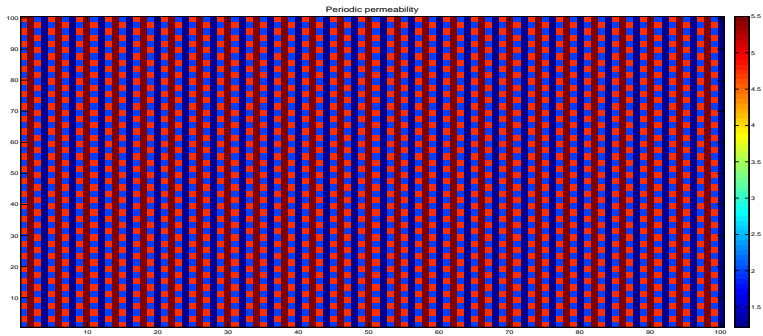}
\caption{Periodic permeability}
\label{exam3-perm}
\end{figure}

\begin{figure}[tbp]
\centering
\includegraphics[width=6in, height=4in]{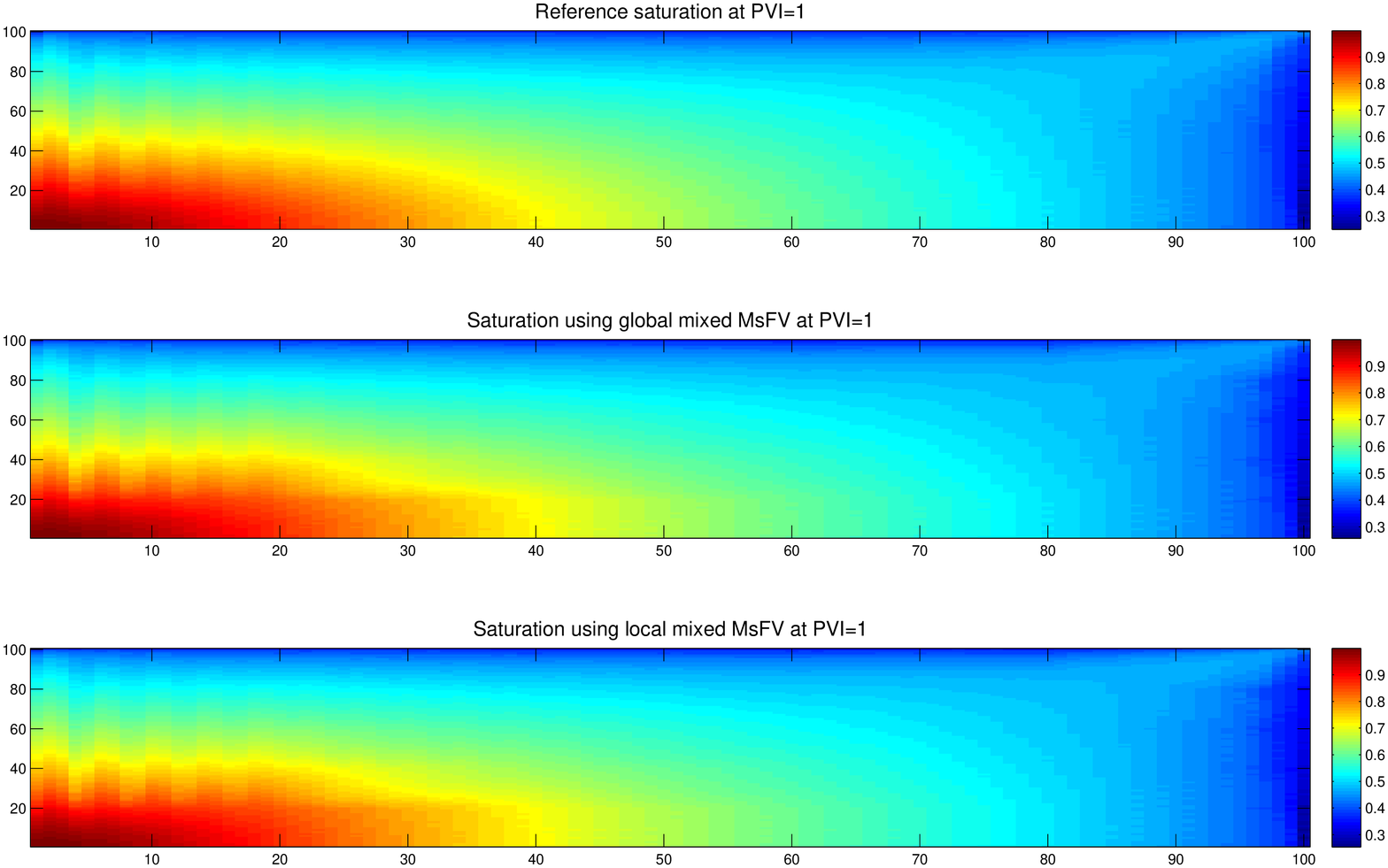}
\caption{Top: Reference saturation at PVI=1 for the third example.
Middle: Saturation at PVI=1 by the global mixed MsFV. Bottom: Saturation at PVI=1 by the local mixed MsFV.}
\label{exam3-saturation}
\end{figure}

\begin{figure}[tbp]
\centering
\includegraphics[width=4in, height=2in]{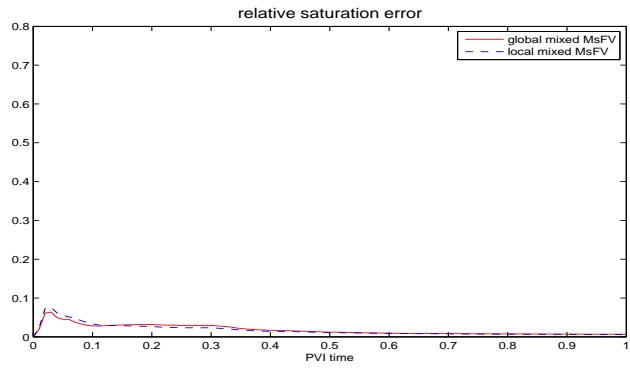}
\caption{saturation error via time for the third example}
\label{exam3-saturation-error}
\end{figure}

\begin{figure}[tbp]
\centering
\includegraphics[width=4in, height=2in]{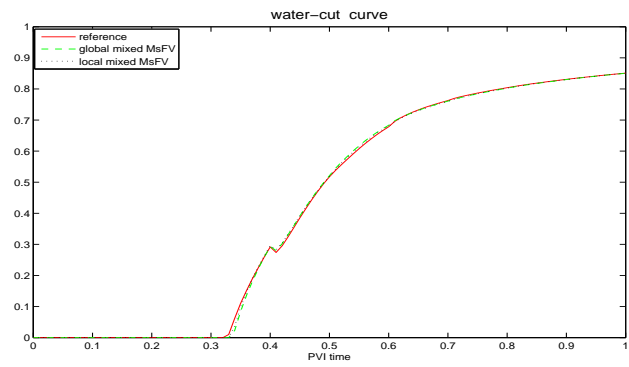}
\caption{water-cut curves for the third example}
\label{exam3-WC}
\end{figure}

In the third    numerical example, we choose the  permeability which is specified by
a periodic function
 \[
 k(x,y)=\frac{2+1.8 \sin(2 \pi x/\epsilon)} {2+1.8 \sin(2\pi y/\epsilon)} + \frac{2+1.8 \sin(2\pi y/\epsilon)}{2+1.8 \cos(2 \pi x /\epsilon)}, \ \  \epsilon=1/25.
 \]
The permeability is defined on $100\times 100$  fine grid and its map is depicted in Figure \ref{exam3-perm}.
Figure \ref{exam3-saturation} depicts the
reference (fine-scale) saturation, the saturation
 using the global
mixed MsFV and the saturation  using the local mixed MsFV, respectively, at PVI=1.
 Here $5\times 5$ coarse grid is taken for simulation in
 both the global mixed MsFV and the local mixed MsFV and the viscosity  ratio
is $\mu_w/\mu_o=1/10$.   Figure  \ref{exam3-saturation} shows that the
saturation profile by global mixed MsFV is almost the same as
the saturation profile by local mixed MsFV and both of them are pretty close to the
reference saturation profile.  The saturation error
via different PVI times  is shown in Figure \ref{exam3-saturation-error}.
It can be seen that the two saturation errors are very close and quite small.
The water-cut curves are shown in Figure \ref{exam3-WC}
for reference, global mixed MsFV and local mixed MsFV, respectively.
Figure \ref{exam3-WC} shows that the three water-cut curves coincides each other.
The average errors of saturation and water-cut
are shown in Table \ref{exam3-tab}.  We find from Table \ref{exam3-tab} that:
(1) the average errors of saturation and water-cut by global mixed MsFV
and local mixed MsFV are comparable. (2) the average errors of saturation are less than
$1\%$ and the average errors of water-cut are less than $2\%$.
This example shows that local mixed MsFV can provide very good accuracy in
the case of separable scales (e.g., periodic) and its performance is as good as the  global
mixed MsFV. The example confirms the findings for many multiscale finite elements (e.g.,\cite{arbogast02, ch03,jym09}).

\begin{table}[hbtp]
\centering \caption{relative errors for saturation and water-cut for the third example}
\begin{tabular}{|c|c|c|c|c|}
\hline
mixed MsFV  &  Water-Cut Error   &  Saturation Error \\

\hline
local mixed MsFV   & 0.0074 & 0.0171\\
\hline
global mixed MsFV  & 0.0082  & 0.0169\\
\hline
\end{tabular}
\label{exam3-tab}
\end{table}

\section{Conclusions}
We developed a framework of   mixed MsFV methods
for elliptic equations arising from flows in porous media.  These
methods take advantages of both mixed MsFEM and finite volume methods.
The essential characteristic of mixed MsFv is the explicit approximation of both pressure and velocity.
We proposed a new multiscale basis
functions for velocity.  Global information can be used
to construct the multiscale velocity basis functions to improve accuracy for highly  heterogeneous porous
media. We  analyze  one of the mixed MsFV  methods and apply it  to simulate two-phase flows in  porous media.
We test the method on permeability fields with both non-separable and separable scales.
Numerical examples demonstrate that the mixed MsFV
can efficiently approximate the two-phase flows on coarse grid. Using global
information in the mixed MsFV yields much better accuracy than the local mixed MsFV if the permeability
field   has strong no-local features.

We also can use multiscale basis for the pressure.
Discrete harmonic and energy minimizing basis constructed in a multilevel fashion is accurate and the
computations are very efficient \cite{XuZikatanov:2004, dmz09}. Moreover,
the methods easily can be extended to unstructured grids without requiring extra geometric information.
 Further investigation of these issues  is worth pursuing in the future.

\appendix\section{Upwind finite volume method for Equation \ref{fsseqn} }
\label{app1}
In the Appendix,  we would like to present  a finite volume
discretization of  the saturation equation (\ref{fsseqn}).  Let $\gamma_{ij}$ be the common face (or edge)  of $K_i^h$ (underlying fine grid) and $K_j^h$ (underlying find grid) and $n_{ij}$ be the normal
vector pointing from $K_i^h$ to $K_j^h$.
Using the $\theta$-rule for temporal discretization and a finite-volume
scheme for  the saturation equation, it follows the following form.
\begin{equation}
\label{numerical-sat}
\frac{1}{\Delta t} (S_i^{n+1}-S_i^n)+\frac{1}{|K_i^h|}\sum_{j\neq i}[\theta F_{ij}(S^{n+1})+(1-\theta)F_{ij}
(S^n)]=0,
\end{equation}
where $S_i^n$ is the cell-average of water saturation at $t=t_n$, i.e.,
\[
S_i^n=\langle S(x,t_n)\rangle_{K_i^h}
\]
and $F_{ij}$ is a numerical approximation of the flux over $\gamma_{ij}$, i.e.,
\[
F_{ij}(S)\approx \int_{\gamma_{ij}} f_w(S)_{ij}u_{ij}\cdot n_{ij} ds.
\]
 Here $f_w(S)_{ij}$ denotes the fractional-flow function associated with $\gamma_{ij}$
and the first-order upstream weighting scheme for it is defined as
\[
f_w(S)_{ij}=\left\{
\begin{array}{ll}
f_w(S_i)  & \mbox{if $u\cdot n_{ij}\geq 0$}\\
f_w(S_j)   & \mbox{if $u\cdot n_{ij}<0$.}
\end{array}
\right.
\]
For $\theta=0$ or  $1$, we can write (\ref{numerical-sat}) as a vector form
\[
S^{n+1}=S^n+(\delta_x^t)^T Af(S^m),  \ \ m=n \ \ or \ \ n+1,
\]
where $(\delta_x^t)_i={\Delta t\over |K_i^h|}$.

If $\theta=0$, then (\ref{numerical-sat}) is an explicit scheme and
only stable provided that time step $\Delta t$ satisfies a stability
condition.

For $\theta=1$, (\ref{numerical-sat}) is an implicit scheme and
unconditionally stable but gives rise to a nonlinear system.
Such a nonlinear system is often solved with a
Newton-Raphson iterative method.
Define
\begin{equation}
\label{Newton-eq1}
G(S^{n+1})=S^{n+1}-S^n-(\delta_x^t)^T Af(S^{n+1}).
\end{equation}
By Taylor expansion, we have
\[
G(S^{n+1})\approx G(S^n)+G'(S^n)(S^{n+1}-S^n).
\]
Noticing $G(S^{n+1})=0$ we have $\delta S^n:=S^{n+1}-S^n=-[G'(S^n)]^{-1} G(S^n)$.
From (\ref{Newton-eq1}), we have
\[
G'(S)=I-(\delta_x^t)^T Af'(S),
\]
where $f'(S)_i=f'(S_i)$.  Hence
\[
S^{n+1}=S^n+\delta S^n.
\]
This iteration proceeds until pre-defined iterations are reached   or  the norm of $\delta S^n$ is smaller than a prescribed value.

\section*{Acknowledgments}
The authors would like to thank the reviewers who provide comments to improve the paper.
L. Jiang  acknowledges  the support from the ExxonMobil Upstream Research Company
for the research.

\end{document}